\documentstyle{amltd2004}
\begin{document}
\annalsline{158}{2003}
\received{May 10, 2000}
\startingpage{419}
\def\bye{\end{document}}
 \font\tenrm=cmr10
\def\ritem#1{\item[{\rm #1}]}
\def\tfrac#1#2{{\textstyle \frac#1#2}}
\input amssym.def
\input amssym.tex
\def\scr#1{{\scriptstyle #1}}
\def\lscr#1{{\scriptstyle #1}\hskip -3pt} 
\def\rscr#1{\hskip-3pt{\scriptstyle #1}}
\def\dnhs{\hskip-8pt}
 \def\iint{\int\hskip-4pt \int}
%--------------- Author macros ---------------
%for Bbb in amstex
\catcode`\@=11
\font\twelvemsb=msbm10 scaled 1100
\font\tenmsb=msbm10
%\font\ninemsb=msbm7 scaled 1100%msbm9
\font\ninemsb=msbm10 scaled 800
\newfam\msbfam
\textfont\msbfam=\twelvemsb  \scriptfont\msbfam=\ninemsb
  \scriptscriptfont\msbfam=\ninemsb
\def\msb@{\hexnumber@\msbfam}
\def\Bbb{\relax\ifmmode\let\next\Bbb@\else
 \def\next{\errmessage{Use \string\Bbb\space only in math
mode}}\fi\next}
\def\Bbb@#1{{\Bbb@@{#1}}}
\def\Bbb@@#1{\fam\msbfam#1}
\catcode`\@=12

 \catcode`\@=11
\font\twelveeuf=eufm10 scaled 1100
\font\teneuf=eufm10
\font\nineeuf=eufm7 scaled 1100%eufm9
\newfam\euffam
\textfont\euffam=\twelveeuf  \scriptfont\euffam=\teneuf
  \scriptscriptfont\euffam=\nineeuf
\def\euf@{\hexnumber@\euffam}
\def\frak{\relax\ifmmode\let\next\frak@\else
 \def\next{\errmessage{Use \string\frak\space only in math
mode}}\fi\next}
\def\frak@#1{{\frak@@{#1}}}
\def\frak@@#1{\fam\euffam#1}
\catcode`\@=12
%-------------- Author entries --------------------

\def\CC{{\Bbb C}}
\def\EE{{\Bbb E}}
\def\FF{{\Bbb F}}
\def\HH{{\Bbb H}}
\def\NN{{\Bbb N}}
\def\QQ{{\Bbb Q}}
\def\RR{{\Bbb R}}
\def\TT{{\Bbb T}}
\def\ZZ{{\Bbb Z}}

\def\H{{\frak H}}
\def\Iso{\operatorname{Iso}}
\def\e{{\rm e}}
\def\i{{\rm i}}
\def\operatorname#1{\mathop{\rm #1}}
\def\cc{\operatorname{cc}}
\def\id{\operatorname{id}}
\def\j{\operatorname{j{}}}
\def\C{\operatorname{C{}}}
\def\G{\operatorname{G{}}}
\def\L{\operatorname{L{}}}
\def\M{\operatorname{M{}}}
\def\GL{\operatorname{GL}}
\def\Op{\operatorname{Op}}
\def\sOp{\widetilde{\operatorname{Op}}}
\def\PGL{\operatorname{PGL}}
\def\Res{\operatorname{Res}}
\def\S{\operatorname{S{}}}
\def\Sp{\operatorname{Sp}}
\def\Sw{{\cal S}}
\def\SL{\operatorname{SL}}
\def\SO{\operatorname{SO}}
\def\PSL{\operatorname{PSL}}
\def\O{\operatorname{O{}}}
\def\T{\operatorname{T{}}}
\def\tr{\operatorname{tr}}
\def\sgn{\operatorname{sgn}}
\def\sign{\operatorname{sign}}
\def\span{\operatorname{span}}
\def\supp{\operatorname{supp}}
\def\Vol{\operatorname{Vol}}
\def\Area{\operatorname{Area}}
\def\ord{\operatorname{ord}}
\def\Prop{\operatorname{Prob}}

\def\I{\operatorname{I}}
\def\II{\operatorname{II}}
\def\III{\operatorname{III}}

\font\mbi=cmmib10 scaled 11
 
\def\veca{{\bf  a}}
\def\vecb{{\bf  b}}
\def\vecc{{\bf  c}}
\def\veck{{\bf  k}}
\def\vecl{{\bf  l}}
\def\vecm{{\bf  m}}
\def\vecn{{\bf  n}}
\def\vecq{{\bf  q}}
\def\vecs{{\bf  s}}
\def\vect{{\bf  t}}
\def\vecu{{\bf  u}}
\def\vecw{{\bf  w}}
\def\vecx{{\bf  x}}
\def\vecy{{\bf  y}}
\def\vecalf{\mathbold{\alpha}}
\def\vecbeta{\mathbold{\beta}}
\def\vecmu{\mathbold{\mu}}
\def\veceta{\mathbold{\eta}}
\def\vecomega{\mathbold{\omega}}
\def\vecsigma{\mathbold{\sigma}}
\def\vecxi{{\mathbold{\xi}}}
\def\veczeta{\mathbold{\zeta}}

\def\vecnull{\mathbold{0}}
\def\vecahalf{{\bf \tfrac12}}

\def\scra{{\cal A}}
\def\scrk{{\cal K}}
\def\scrm{{\cal M}}
\def\scrn{{\cal N}}
\def\scrq{{\cal Q}}
\def\scrs{{\cal S}}
\def\scrt{{\cal T}}
\def\scru{{\cal U}}
\def\scrx{{\cal X}}
\def\scry{{\cal Y}}

\def\Re{\operatorname{Re}\nolimits}
\def\Im{\operatorname{Im}\nolimits}

\def\trans{\,^{\rm t}\!}
%%%%%%%%%%%%%%%%%%%%%%%%%%%%%%%%%%%%%%%%%%%%%%%%%%%%%%%%%%%%%%%%%%%%%%%%%%

\title{Pair correlation densities of\\ inhomogeneous quadratic forms} 
\shorttitle{Inhomogeneous quadratic forms}  

 \author{Jens Marklof}
 \institutions{School of Mathematics, University of Bristol,
Bristol, United Kingdom\\
{\eightpoint {\it E-mail address\/}: j.marklof@bristol.ac.uk
}}

\centerline{\bf Abstract}
\vglue8pt

Under explicit diophantine conditions on $(\alpha,\beta)\in\RR^2$, we 
prove that the local two-point correlations of the
sequence given by the values $(m-\alpha)^2+\break (n-\beta)^2$, with 
$(m,n)\in\ZZ^2$, are those of a Poisson process.
This partly confirms a conjecture
of Berry and Tabor [2] on spectral statistics of quantized
integrable systems, and also establishes a particular case of the 
quantitative version of the Oppenheim conjecture for
inhomogeneous quadratic forms of signature (2,2).
The proof uses theta sums and Ratner's classification of measures
invariant under unipotent flows.
 
\vglue-6pt
\section{Introduction}

 \vglue-6pt
1.1.
Let us denote by $0\leq\lambda_1\leq\lambda_2\leq \cdots\rightarrow\infty$
the infinite sequence given by the values of 
$$
(m-\alpha)^2 + (n-\beta)^2
$$
at lattice points $(m,n)\in\ZZ^2$, for fixed $\alpha,\beta\in[0,1]$.
In a numerical experiment, Cheng and Lebowitz \cite{Cheng91}
found that, for generic $\alpha,\beta$, 
the local statistical measures of the {\it deterministic} sequence $\lambda_j$
appear to be those of independent {\it random} variables from a
Poisson process. 

\vglue6pt 1.2.  
This numerical observation supports
a conjecture of Berry and Tabor \cite{Berry77} in the context of quantum chaos,
according to which the local eigenvalue statistics of generic quantized
integrable systems are Poissonian. In the case discussed here,
the $\lambda_j$ may be viewed
(up to a factor $4\pi^2$) as the eigenvalues of the Laplacian
$$
-\Delta = -\frac{\partial^2}{\partial x^2}-\frac{\partial^2}{\partial y^2}
$$
with quasi-periodicity conditions
$$
\varphi(x+k,y+l) = \e^{-2\pi\i(\alpha k+\beta l)} \varphi(x,y) , \quad
k,l\in\ZZ.
$$
The corresponding classical dynamical system is the geodesic flow on the 
unit tangent bundle of the flat torus \pagebreak $\TT^2$.

 1.3.
The asymptotic density of the sequence of $\lambda_j$ is $\pi$, according
to the well known formula for the number of lattice points
in a large, shifted circle:
$$
\#\{ j : \lambda_j\leq\lambda \}
=\#\{ (m,n)\in\ZZ^2: (m-\alpha)^2 + (n-\beta)^2 \leq \lambda \}
\sim \pi \lambda  
$$
for $\lambda\rightarrow\infty$. The rate of convergence is discussed
in detail by Kendall \cite{Kendall48}.

\vglue12pt 1.4.
More generally, suppose we have a sequence 
$\lambda_1\leq\lambda_2\leq\cdots\rightarrow\infty$ of mean density $D$, i.e.,
$$
\lim_{\lambda\rightarrow\infty}
\frac1\lambda\#\{ j : \lambda_j\leq\lambda \} = D .
$$
For a given interval $[a,b]\subset\RR$, the 
{\it pair correlation function} is then defined as 
$$
R_2[a,b](\lambda) = \frac{1}{D\lambda}  \#\{ j\neq k : \lambda_j\leq\lambda,\;
\lambda_k\leq\lambda,\;
a\leq \lambda_j-\lambda_k \leq b \} .
$$

The following result is classical.

\nonumproclaim{1.5. Theorem} \hskip-9pt  If the $\lambda_j$ come from a Poisson process with
 mean density~$D${\rm ,} 
$$
\lim_{\lambda\rightarrow\infty} R_2[a,b](\lambda)  = D (b-a) 
$$
almost surely.
\endproclaim

1.6.
We will assume throughout most of the paper that $\alpha,\beta,1$ 
are linearly independent over $\QQ$. This makes sure that there are no
systematic degeneracies in the sequence, which would contradict
the independence we wish to establish. The symmetries leading
to those degeneracies can, however, be removed without much difficulty. 
This will be illustrated in Appendix~A.

\vglue12pt 1.7.
We shall need a mild diophantine condition on $\alpha$.
An irrational number $\alpha\in\RR$ is called {\it diophantine}
if there exist constants $\kappa, C>0$ such that
$$
\big| \alpha - \frac{p}{q} \big| > \frac{C}{q^\kappa}
$$
for all $p,q\in\ZZ$. The smallest possible value of $\kappa$
is $\kappa=2$ \cite{Schmidt72}. We will say $\alpha$ is of {\it type $\kappa$}.

\nonumproclaim{1.8. Theorem} Suppose $\alpha,\beta,1$ are linearly independent over $\QQ${\rm ,} and
assume 
$\alpha$ is diophantine.
Then
$$
\lim_{\lambda\rightarrow\infty} R_2[a,b](\lambda)  = \pi (b-a) .
$$
\endproclaim

This proves the Berry-Tabor conjecture for the spectral two-point
correlations of the Laplacian in 1.2.

It is well known that
almost all $\alpha$ (in the measure-theoretic sense) 
are diophantine \cite{Schmidt72}. We therefore have the following corollary.

\nonumproclaim{1.9. {C}orollary} Let $\alpha,\beta$ be independent uniformly distributed random
variables  in $[0,1]$. Then
$$
\lim_{\lambda\rightarrow\infty} R_2[a,b](\lambda) = \pi (b-a) 
$$
almost surely.
\endproclaim

1.10. {\it Remark}.
In \cite{Cheng94}, Cheng, Lebowitz and Major proved convergence
of the expectation value\footnote{They consider a slightly
different statistic, the number of lattice points in
a random circular strip of fixed area. The variance of this
distribution is very closely related to our pair correlation function.} 
$$
\lim_{\lambda\rightarrow\infty} \EE R_2[a,b](\lambda) = \pi (b-a),
$$
that is, on average over $\alpha,\beta$.

\demo{{\rm 1.11.} Remark}
Notice that Theorem 1.8 is much stronger than
the corollary. It provides explicit examples of ``random''
deterministic sequences that satisfy the pair correlation conjecture.
An admissible choice is for instance 
$\alpha=\sqrt2$, $\beta=\sqrt3$ \cite{Schmidt72}.
\enddemo

1.12.
The statement of Theorem 1.8 
does not hold for any rational $\alpha,\beta$,
where the pair correlation function is unbounded 
(see Appendix~A.10 for details). 
This can be used to show that for generic $(\alpha,\beta)$
(in the topological sense) the pair correlation function 
does not converge to a uniform density:

\nonumproclaim{1.13. Theorem} For any $a>0${\rm ,} there exists a set $C\subset\TT^2$ of second Baire
category{\rm ,} for which the following holds.\footnote{A set of
first Baire category is a countable union of nowhere dense sets.
Sets of second category are all those sets  which are not
of first category.}

{\rm (i)} For $(\alpha,\beta)\in C$, there exist arbitrarily large $\lambda$
such that
$$
R_2[-a,a](\lambda)  \geq \frac{\log\lambda}{\log\log\log\lambda}.
$$

{\rm (ii)} For $(\alpha,\beta)\in C${\rm ,} there exists an infinite  sequence
$L_1<L_2<\cdots\rightarrow\infty$ such that
$$
\lim_{j\rightarrow\infty} R_2[-a,a](L_j)  = 2\pi a .
$$
\endproclaim

In the above, $\log\log\log\lambda$ may be replaced by any
slowly increasing positive function $\nu(\lambda)\leq \log\log\log\lambda$
with $\nu(\lambda)\rightarrow\infty$ ($\lambda\rightarrow\infty$).

\vglue4pt 1.14.
The above results can be extended to the pair correlation densities
of forms $(m_1-\alpha_1)^2+\ldots+(m_k-\alpha_k)^2$ in more than
two variables; see \cite{Marklof00} for details.

\pagebreak 1.15. {\it A brief review}.
After its formulation in 1977, Sarnak \cite{Sarnak97} was the first to prove
the Berry-Tabor conjecture for the pair correlation of almost all 
positive definite binary
quadratic forms
$$
\alpha m^2 + \beta mn + \gamma n^2, \quad m,n\in\ZZ 
$$
(``almost all'' in the measure-theoretic sense).
These values represent the eigenvalues of the Laplacian on a flat torus. 
His proof uses averaging
techniques to reduce the pair correlation problem to
estimating the number of solutions of systems of diophantine equations.
The almost-everywhere result then follows from a variant of the
Borel-Cantelli argument. 
For further related examples of sequences whose pair correlation function
converges to the uniform density almost everywhere in parameter space, see
\cite{Rudnick98}, \cite{Rudnick99}, \cite{VanderKam96}, \cite{VanderKam97}, \cite{Zelditch98}. Results
on higher correlations have been obtained 
recently in \cite{Rudnick99a}, \cite{Rudnick99b}, \cite{VanderKam00}.

Eskin, Margulis and Mozes \cite{Eskin98b}
have recently given explicit diophantine conditions 
under which the pair correlation
function of the above binary\break quadratic forms is Poisson. 
Their approach uses ergodic-theoretic methods based on Ratner's
classification of measures invariant under unipotent flows.
This will also be the key ingredient in our proof
for the inhomogeneous set-up. New in the approach
presented here is the application of theta sums 
\cite{Marklof96}, \cite{Marklof96b}, \cite{Marklof98}.

The pair correlation problem for binary quadratic forms
may be viewed  as a special case of the quantitative version of the Oppenheim
conjecture for forms of signature (2,2), which is particularly difficult
\cite{Eskin98}.

\vglue6pt {\it Acknowledgments}.
I thank A.~Eskin, F.~G\"otze, G.~Margulis, S.~Mo\-zes, Z.~Rudnick
and N.~Shah for very helpful discussions and correspondence.
Part of this research was carried out during visits at
the Universities of Bielefeld and Tel Aviv, with financial
support from SFB 343 ``Diskrete Strukturen in der Mathematik''
and the Hermann Minkowski Center for Geometry, respectively.
I have also highly appreciated the referees'  and A.~Str\"ombergsson's 
comments and suggestions on the first version of this paper.

\vglue-18pt
\section{The plan}\label{secplan} 
\vglue-8pt

2.1.
The plan is first to  smooth the pair correlation function,
i.e., to consider
$$
R_2(\psi_1,\psi_2,h,\lambda)
=\frac{1}{\pi\lambda} \sum_{j,k} 
\psi_1\left(\frac{\lambda_j}{\lambda}\right)\psi_2\left(\frac{\lambda_k}{\lambda}\right)
\hat h(\lambda_j-\lambda_k) .
$$
Here $\psi_1,\psi_2\in\Sw(\RR_+)$ are 
real-valued, and $\Sw(\RR_+)$ denotes the Schwartz
class of infinitely differentiable functions
of the half line $\RR_+$ (including the origin), which, as well as their
derivatives, decrease rapidly at $+\infty$. 
It is helpful to think of $\psi_1,\psi_2$ as smoothed characteristic
functions, i.e., positive and with compact support. Note that \pagebreak
$\hat h$ is the Fourier transform of a compactly supported function
$h\in\C(\RR)$, defined by
$$
\hat h(s)= \int_\RR h(u) e(\tfrac12 us)\, du,
$$
with the shorthand $e(z):=\e^{2\pi\i z}$.

We will prove the following (Section \ref{secmain}).

\nonumproclaim{2.2. Theorem}  Let $\psi_1,\psi_2\in\Sw(\RR_+)$ be real\/{\rm -}\/valued{\rm ,}
and $h\in\C(\RR)$ with compact support.
Suppose $\alpha,\beta,1$ are linearly independent over $\QQ${\rm ,} and assume 
$\alpha$ is diophantine.
Then
$$
\lim_{\lambda\rightarrow\infty}
R_2(\psi_1,\psi_2,h,\lambda) 
= \big\{ \hat h(0) + \pi \int_\RR \hat h(s)\,ds \big\}
\int_0^\infty \psi_1(r)\psi_2(r)\, dr.
$$
\endproclaim

The first term comes straight from the terms $j=k$; the second
one is the more interesting.

Theorem~2.2 implies Theorem 1.8 by a 
standard approximation argument (Section \ref{secmain}).

\vglue12pt 2.3. 
Using the Fourier transform we may write
$$
R_2(\psi_1,\psi_2,h,\lambda)
=\frac{1}{\pi\lambda} \int_\RR 
\big(\sum_{j} \psi_1\big(\frac{\lambda_j}{\lambda}\big) e(\tfrac12 \lambda_j u) \big)
\overline{
\big(\sum_{j} \psi_2\big(\frac{\lambda_j}{\lambda}\big) e(\tfrac12 \lambda_j u) \big)}
h(u)\, du .
$$
We will show that the inner sums can be viewed as a theta sum
(see 4.14 for details)
$$
\theta_\psi(u,\lambda)
=\frac{1}{\sqrt\lambda}
\sum_{j} \psi\left(\frac{\lambda_j}{\lambda}\right) e(\tfrac12 \lambda_j u)
$$
living on a certain manifold $\Sigma$ of finite volume 
(Sections \ref{secschrodinger} and 4). 
The integration in
$$
R_2(\psi_1,\psi_2,h,\lambda) 
= \frac1\pi \int_\RR \theta_{\psi_1}(u,\lambda)
\overline{\theta_{\psi_2}(u,\lambda)} h(u)\, du
$$
will then be identified with an orbit of a unipotent flow on
$\Sigma$, which becomes equidistributed as $\lambda\rightarrow\infty$.
The equidistribution follows from Ratner's classification of
measures invariant under the unipotent flow (Section \ref{secunipotent}). 
A crucial subtlety
is that $\Sigma$ is noncompact, and that the theta sum is unbounded
on this noncompact space. This requires careful estimates which
guarantee that no positive mass of the above integral over a small arc
of the orbit escapes to infinity (Section \ref{secdiophantine}). 

The only exception is a small neighbourhood of $u=0$, where in fact 
a positive mass escapes to infinity, giving a contribution
$$
2\pi^2 h(0) \int_0^\infty \psi_1(r) \psi_2(r)\, dr 
= \pi^2 \int_\RR \hat h(s)\,ds \int_0^\infty \psi_1(r) \psi_2(r)\, dr ,
$$
which is the second term in Theorem~2.2.

The remaining part of the orbit becomes equidistributed under
the above diophantine conditions, which yields
$$
\frac{1}{\mu(\Sigma)}\int_\Sigma \theta_{\psi_1}
\overline{\theta_{\psi_2}} d\mu \; \int_\RR h(u)\, du ,
$$
where $\mu$ is the invariant measure (Section \ref{secequi}).
The first integral can be calculated quite easily (Section \ref{secmain}). 
It is
$$
\frac{1}{\mu(\Sigma)}\int_\Sigma  \theta_{\psi_1}
\overline{\theta_{\psi_2}} d\mu \; \int_\RR h(u)\, du 
= \pi \int_0^\infty \psi_1(r) \psi_2(r)\, dr \; \int_\RR h(u)\, du ,
$$
which finally yields
$$
\pi \hat h(0) \int_0^\infty \psi_1(r) \psi_2(r)\, dr ,
$$
the first term in Theorem~2.2.

The proof of Theorem~1.13, which provides a set of counterexamples
to the convergence to uniform density,
is given in Section \ref{counter}.

\section{Schr\"odinger and Shale-Weil representation}\label{secschrodinger}

3.1. Let $\omega$ be the standard symplectic form on $\RR^{2k}$, i.e., 
$$
\omega(\vecxi,\vecxi')=\vecx\cdot\vecy'-\vecy\cdot\vecx' ,
$$
where 
$$
\vecxi=
\left( \begin{array}{c}
\vecx \\
\vecy \\
\end{array}\right), 
\qquad
\vecxi'=
\left( \begin{array}{c}
\vecx' \\
\vecy' \\
\end{array}\right), \qquad
\vecx,\vecy,\vecx',\vecy'\in\RR^k .
$$

The {\it Heisenberg group} $\HH(\RR^k)$ is then defined as the set 
$\RR^{2k}\times\RR$ with multiplication law \cite{Lion80} 
$$
(\vecxi, t)(\vecxi', t')=(\vecxi+\vecxi', 
t+t'+\tfrac12\omega(\vecxi,\vecxi') ) .
$$
Note that we have the decomposition
$$
\left(\left( \begin{array}{c}
\vecx \\
\vecy \\
\end{array}\right), t\right)
=
\left(\left( \begin{array}{c}
\vecx \\
\vecnull \\
\end{array}\right), 0\right)
\left(\left( \begin{array}{c}
\vecnull \\
\vecy \\
\end{array}\right), 0\right)
(\vecnull , t-\tfrac12 \vecx\cdot\vecy).
$$

3.2.  
The {\it Schr{\rm \"{\it o}}dinger representation} of $\HH(\RR^k)$
on $f\in\L^2(\RR^k)$
is given by (cf.~\cite[p.\ 15]{Lion80})
$$\begin{array}{rll}
\left[W\left(\left( \begin{array}{c}
\vecx \\
\vecnull \\
\end{array}\right), 0\right) f\right](\vecw) &= e(\vecx\cdot\vecw)\; f(\vecw),
&\quad \hbox{with $\vecx,\vecw\in\RR^k$,}
\\ \noalign{\vskip4pt}
 \left[W\left(\left( \begin{array}{c}
\vecnull \\
\vecy \\
\end{array}
\right), 0\right) f \right] (\vecw)& =  f(\vecw-\vecy),
&\quad \hbox{with $\vecy,\vecw\in\RR^k$,}
\\  \noalign{\vskip4pt}
W(\vecnull , t)& = e(t) \id,
&\quad \hbox{with $t\in\RR$.} \end{array}
$$
Therefore for a general element $(\vecxi,t)$ in $\HH(\RR^k)$ 
$$
\left[W\left(\left( \begin{array}{c}
\vecx \\
\vecy \\
\end{array}\right), t\right) f\right](\vecw) = e(t-\tfrac12\vecx\cdot\vecy)
\; e(\vecx\cdot\vecw)\; f(\vecw-\vecy) .
$$

\vglue12pt 3.3. 
For every element $M$ in the symplectic group $\Sp(k,\RR)$ of $\RR^{2k}$,
we can define a new representation $W_M$ of $\HH(\RR^k)$ by
$$
W_M(\vecxi,t)=W(M\vecxi,t) .
$$
All such representations are irreducible and,
by the Stone-von Neumann theorem, unitarily equivalent (see \cite{Lion80}
for details). That is, for each $M\in\Sp(k,\RR)$ there exists 
a unitary operator $R(M)$ such that
$$
R(M)\; W(\vecxi,t)\; R(M)^{-1} = W(M\vecxi,t) .
$$
The $R(M)$ is determined up to a unitary phase factor and defines
the projective {\it Shale-Weil representation} of the symplectic group.
{\it Projective} means that
$$
R(MM')=c(M,M')R(M)R(M')
$$
with cocycle 
$c(M,M')\in\CC$, $|c(M,M')|=1$, but $c(M,M')\neq 1$ in general.

\vglue12pt 3.4. 
For our present purpose it suffices to consider the group
$\SL(2,\RR)$ which is embedded in $\Sp(k,\RR)$ by 
$$
\left( \begin{array}{cc}
a & b \\
c & d
\end{array}\right)
\mapsto
\left( \begin{array}{cc}
a\, 1_k & b\, 1_k \\
c\, 1_k & d\, 1_k
\end{array}\right)
$$ 
where $1_k$ is the $k\times k$ unit matrix.

The action of
$M\in\SL(2,\RR)$ on $\vecxi\in\RR^{2k}$ is then given by
$$
M\vecxi = \left( \begin{array}{c}
a\vecx +b\vecy \\
c\vecx +d\vecy \\
\end{array}\right), \quad \hbox{with }
M=\left( \begin{array}{cc}
a&b\\
c&d
\end{array}\right),\quad 
\vecxi= \left( \begin{array}{c}
\vecx\\
\vecy
\end{array}\right). 
$$
 \vglue12pt
 3.5. For $M\in\SL(2,\RR)\hookrightarrow\Sp(k,\RR)$ 
we have the explicit representations 
(see \cite[p.\ 61f]{Lion80}).
\begin{eqnarray*}
&&\hskip-.5in
[R(M)f](\vecw)\\[4pt]
=
&&\left\{ \begin{array}{ll}
\displaystyle
|a|^{k/2} e(\tfrac12 \|\vecw\|^2 ab) f(a\vecw) & (c=0)\\[2mm]
\displaystyle
|c|^{-k/2} \int_{\RR^k} e\left[\frac{\tfrac12 (a\|\vecw\|^2+d\|\vecw'\|^2)-
\vecw\cdot\vecw'}{c}\right]\, f(\vecw')\,d\vecw' & (c\neq 0).
\end{array}\right. 
\end{eqnarray*}
Here $\|\;\cdot\;\|$ denotes the euclidean norm in $\RR^k$,
$$
\|\vecx\|=\sqrt{x_1^2+\cdots+x_k^2}. \pagebreak
$$

 3.6.  
If
$$
M_1=\left( \begin{array}{cc}
a_1&b_1\\
c_1&d_1
\end{array}\right),
\quad
M_2=\left( \begin{array}{cc}
a_2&b_2\\
c_2&d_2
\end{array}\right),
\quad
M_3=\left( \begin{array}{cc}
a_3&b_3\\
c_3&d_3
\end{array}\right) \in\SL(2,\RR),
$$
with $M_1M_2=M_3$,
the corresponding cocycle is
$$
c(M_1,M_2)=\e^{-\i\pi k\sign(c_1c_2c_3)/4} ,
$$
where
$$
\sign(x)=
\left\{ \begin{array}{ll}
-1 & (x<0)\\
0 & (x=0)\\
1& (x>0).
\end{array}\right. 
$$

\vglue12pt 3.7.  
In the special case when
$$
M_1=\left( \begin{array}{cc}
\cos\phi_1 & -\sin\phi_1\\
\sin\phi_1 &  \cos\phi_1
\end{array}\right),
\quad
M_2=\left( \begin{array}{cc}
\cos\phi_2 & -\sin\phi_2\\
\sin\phi_2 &  \cos\phi_2
\end{array}\right),
$$
we find
$$
c(M_1,M_2)=\e^{-\i\pi k(\sigma_{\phi_1}+\sigma_{\phi_2}
-\sigma_{\phi_1+\phi_2})/4}
$$
where
$$
\sigma_\phi=
\left\{ \begin{array}{ll}
2\nu   & \hbox{if $\phi=\nu\pi$,}\\
2\nu+1 & \hbox{if $\nu\pi<\phi<(\nu+1)\pi$.}
\end{array}\right. 
$$

\vglue12pt 3.8.  Every $M\in\SL(2,\RR)$ admits the unique Iwasawa decomposition
$$
M = \left( \begin{array}{cc}
1 & u\\
0 & 1
\end{array}\right)
\left( \begin{array}{cc}
v^{1/2} & 0\\
0 & v^{-1/2}
\end{array}\right)
\left( \begin{array}{cc}
\cos\phi & -\sin\phi\\
\sin\phi &  \cos\phi
\end{array}\right)
=
(\tau,\phi),
$$
where $\tau=u+\i v\in\H$, $\phi\in[0,2\pi)$.
This parametrization leads to the well known action of
$\SL(2,\RR)$ on $\H\times[0,2\pi)$,
$$
\left( \begin{array}{cc}
a & b \\
c & d
\end{array}\right)
(\tau,\phi)= (\frac{a\tau+b}{c\tau+d}, \phi+\arg(c\tau+d) \bmod2\pi) .
$$
We will sometimes use the convenient notation 
$(M\tau,\phi_M):=M(\tau,\phi)$ and
$u_M:=\Re(M\tau)$, $v_M:=\Im(M\tau)$.

\vglue12pt 3.9. 
The (projective) Shale-Weil representation of $\SL(2,\RR)$ reads 
in these coordinates
$$
[R(\tau,\phi)f](\vecw)=[R(\tau,0)R(\i,\phi)f](\vecw)
=
v^{k/4} e(\tfrac12 \|\vecw\|^2 u) [R(\i,\phi)f](v^{1/2}\vecw)
$$
and
\begin{eqnarray*}
&\hskip-20pt&[R(\i,\phi)f](\vecw)\\[4pt]
&\hskip-20pt&\qquad =
\left\{ \begin{array}{ll}
f(\vecw) & (\phi=0\bmod 2\pi)\\[2mm]
f(-\vecw) & (\phi=\pi\bmod 2\pi)\\[2mm]
\displaystyle
|\sin\phi|^{-k/2} 
\int_{\RR^k} e\left[\frac{\tfrac12(\|\vecw\|^2+\|\vecw'\|^2)\cos\phi-
\vecw\cdot\vecw'}{\sin\phi}\right]& f(\vecw')\,d\vecw'  \\
& (\phi\neq 0\bmod\pi).
\end{array}\right. 
\end{eqnarray*}
Note that $R(\i,\pi/2)={\cal F}$ is the Fourier transform.

\vglue12pt 3.10. 
For Schwartz functions $f\in\Sw(\RR^k)$,
\begin{eqnarray*}
&&
\lim_{\phi\rightarrow 0\pm} |\sin\phi|^{-k/2} 
\int_{\RR^k} e\left[\frac{\tfrac12(\|\vecw\|^2+\|\vecw'\|^2)\cos\phi-
\vecw\cdot\vecw'}{\sin\phi}\right]\, f(\vecw')\,d\vecw' 
\\[5pt]
&&\qquad\qquad =\e^{\pm\i\pi k \pi/4} f(\vecw),
\end{eqnarray*}
and hence
this projective representation is in general discontinuous at $\phi=\nu\pi$,
$\nu\in\ZZ$. This can be overcome by setting
$$
\tilde R(\tau,\phi) = \e^{-\i\pi k\sigma_\phi/4} R(\tau,\phi) .
$$
In fact, $\tilde R$ corresponds  to a unitary representation
of the double cover of $\SL(2,\RR)$ \cite{Lion80}. This means in particular
that (compare 3.7)
$$
\tilde R(\i,\phi)\tilde R(\i,\phi')=\tilde R(\i,\phi+\phi'),
$$
where $\phi\in[0,4\pi)$ parametrizes the double cover of 
$\SO(2)\subset\SL(2,\RR)$.

 \section{Theta sums\label{sectheta}}

4.1. 
The {\it Jacobi group} is defined as the semidirect product 
\cite{Berndt98}
$$
\Sp(k,\RR)\ltimes\HH(\RR^k)
$$
with multiplication law
$$
(M;\vecxi,t)(M';\vecxi',t')=
(MM';\vecxi+M\vecxi', t+t'+\tfrac12\omega(\vecxi,M\vecxi')) .
$$
This definition is motivated by the fact that, since 
$$
R(M) W(\vecxi',t') =W(M\vecxi',t')R(M),
$$
(recall 3.3) we have
\begin{eqnarray*}
&&\hskip-.5in W(\vecxi,t)R(M)\;W(\vecxi',t')R(M')\\
&=&W(\vecxi,t)W(M\vecxi',t')\;R(M)R(M') \\
&=&c(M,M')^{-1}
\;W(\vecxi+M\vecxi',t+t'+\tfrac12\omega(\vecxi,M\vecxi'))\; R(MM').
\end{eqnarray*}
Hence
$$
R(M;\vecxi,t)=W(\vecxi,t)R(M)
$$
defines a projective representation of the Jacobi group,
with cocycle $c(M,M')$ as above,
the so-called {\it Schr{\rm \"{\it o}}dinger\/{\rm -}\/Weil representation}
\cite{Berndt98}.

Let us also put
$$
\tilde R(\tau,\phi;\vecxi,t)=W(\vecxi,t)\tilde R(\tau,\phi).
$$

\vglue-9pt
\phantom{1}

\demo{{\rm 4.2.} Jacobi\/{\rm '}\/s theta sum}
We define Jacobi's theta sum for $f\in{\cal S}(\RR^k)$ by
$$
\Theta_f(\tau,\phi;\vecxi,t)= \sum_{\vecm\in\ZZ^k}
[\tilde R(\tau,\phi;\vecxi,t) f](\vecm) .
$$
More explicitly, for $\tau=u+\i v$, $\vecxi =\left( \begin{array}{c}
\vecx \\
\vecy \\
\end{array}\right)$,
$$
\Theta_f(\tau,\phi;\vecxi,t) = 
v^{k/4} e(t-\tfrac12 \vecx\cdot\vecy)
\sum_{\vecm\in\ZZ^k} f_\phi( (\vecm-\vecy) v^{1/2}) 
e(\tfrac12\|\vecm-\vecy\|^2 u + \vecm\cdot\vecx) ,
$$
where
$$
f_\phi = \tilde R(\i,\phi)f .
$$
It is easily seen that if $f\in{\cal S}(\RR^k)$ then
$f_\phi\in{\cal S}(\RR^k)$ for $\phi$ fixed,
and thus also $\tilde R(\tau,\phi;\vecxi,t) f\in{\cal S}(\RR^k)$
for fixed $(\tau,\phi;\vecxi,t)$.
This guarantees rapid convergence of the above series.
We have the following uniform bound.
\enddemo

\nonumproclaim{4.3. Lemma}  Let $f_\phi=\tilde R(\i,\phi) f${\rm ,} with $f\in{\cal S}(\RR^k)$. Then{\rm ,}
for any $R>1${\rm ,} there is a constant $c_R$ such that for all $\vecw\in\RR^k${\rm ,}
$\phi\in\RR${\rm ,} 
$$
|f_\phi(\vecw)| \leq c_R (1+\|\vecw\|)^{-R} .
$$
\endproclaim

\demo{Proof}
Since $f\in{\cal S}(\RR^k)$, we can use repeated integration by parts to
show that 
$$
\left| |\sin\phi|^{-k/2} 
\int_{\RR^k} e\left[\frac{\tfrac12(\|\vecw\|^2+\|\vecw'\|^2)\cos\phi-
\vecw\cdot\vecw'}{\sin\phi}\right]  f(\vecw')\,d\vecw' \right| 
\leq 
c'_R (1+\|\vecw\|)^{-R} 
$$
uniformly for all $\phi\notin (\nu\pi-\frac{1}{100},\nu\pi+\frac{1}{100})$,
$\nu\in\ZZ$. That is, 
$$
|f_\phi(\vecw)| \leq c'_R (1+\|\vecw\|)^{-R} 
$$
in the above range.

Furthermore $f_{\pi/2}$ is up to a phase factor $\e^{\i\pi k}$
the Fourier transform of $f$ and therefore of Schwartz class as well.
Again, after integration by parts,
\begin{eqnarray*}
&&\left| |\sin\phi|^{-k/2} 
\int_{\RR^k} e\left[\frac{\tfrac12(\|\vecw\|^2+\|\vecw'\|^2)\cos\phi-
\vecw\cdot\vecw'}{\sin\phi}\right]\, f_{\pi/2}(\vecw')\,d\vecw' \right| 
\\[4pt]
&&\qquad\qquad \leq 
c''_R (1+\|\vecw\|)^{-R} 
\end{eqnarray*}
for all $\phi\notin (\nu\pi-\frac{1}{100},\nu\pi+\frac{1}{100})$,
$\nu\in\ZZ$.
This means
$$
|f_{\phi+\pi/2}(\vecw)| \leq c''_R (1+\|\vecw\|)^{-R} 
$$
in the above range, or, by replacement of  $\phi\mapsto\phi-\pi/2$,
$$
|f_{\phi}(\vecw)| \leq c''_R (1+\|\vecw\|)^{-R} ,
$$
for all $\phi\notin (\nu\pi+\tfrac12\pi-\frac{1}{100},
\nu\pi+\tfrac12 \pi+\frac{1}{100})$,
$\nu\in\ZZ$.

Clearly for each $\phi\in\RR$ at least one of the bounds applies;
we put $c_R=\max\{c'_R,c''_R\}$.
\enddemo

4.4. The following transformation formulas are crucial for our further
investigations:

\nonumproclaim{Jacobi 1} 
$$
\Theta_f\left(-\frac{1}{\tau},\phi+\arg\tau;
\left( \begin{array}{c}-\vecy \\\vecx\end{array}\right),t\right)
= \e^{-\i\pi k/4} \Theta_f\left(\tau,\phi;
\left( \begin{array}{c}\vecx \\\vecy\end{array}\right),t\right).
$$
\endproclaim

\demo{Proof}
The Poisson summation formula states that for any $f\in{\cal S}(\RR^k)$, 
$$
\sum_{\vecm\in\ZZ^k}
[{\cal F} f](\vecm)
=\sum_{\vecm\in\ZZ^k} f(\vecm)
$$
where ${\cal F}$ is the Fourier transform.
Because
$$
{\cal F}
=R(\i,\pi/2)= R(S), \quad 
S=\left( \begin{array}{cc} 0 & -1  \\ 1 & 0 \end{array}\right),
$$
and secondly $\tilde R(\tau,\phi;\vecxi,t) f\in{\cal S}(\RR^k)$
for fixed $(\tau,\phi;\vecxi,t)$,
the Poisson summation formula yields
$$
\sum_{\vecm\in\ZZ^k}
[R(S)\tilde R(\tau,\phi;\vecxi,t) f](\vecm)
=\sum_{\vecm\in\ZZ^k}[\tilde R(\tau,\phi;\vecxi,t) f](\vecm).
$$
We have
$$
R(S) \tilde R(\tau,\phi;\vecxi,t) 
= R(S) W(\vecxi,t) \tilde R(\tau,0) \tilde R(\i,\phi)
= W(S\vecxi,t) R(S)  R(\tau,0) \tilde R(\i,\phi) ;
$$
furthermore
$$
R(S)  R(\tau,0) = R\left(-\frac1\tau,\arg\tau\right) = R\left(-\frac1\tau,0\right) R(\i,\arg\tau),
$$
since $(\tau,0)$ and $(-\frac1\tau,0)$ are upper triangular matrices, 
and hence the corresponding cocycles are trivial, i.e., 
equal to 1 (recall 3.6).
Finally, since $0<\arg\tau<\pi$ for $\tau\in\H$, 
$$
R(\i,\arg\tau)\tilde R(\i,\phi)
=\e^{\i\pi k/4}\tilde R(\i,\arg\tau)\tilde R(\i,\phi)
=\e^{\i\pi k/4}\tilde R(\i,\phi+\arg\tau).
$$
Collecting all terms, we find
$$
R(S)\tilde R(\tau,\phi;\vecxi,t)
=\e^{\i\pi k/4}\tilde R\left(-\frac1\tau,\phi+\arg\tau;S\vecxi,t\right),
$$
and hence
$$
\sum_{\vecm\in\ZZ^k}
\left[\tilde R\left(-\frac1\tau,\phi+\arg\tau;S\vecxi,t\right) f\right](\vecm)
=\e^{-\i\pi k/4}\sum_{\vecm\in\ZZ^k}[\tilde R(\tau,\phi;\vecxi,t) f](\vecm),
$$
which proves the claim.
\enddemo

\nonumproclaim{Jacobi 2} 
$$
\Theta_f\left(\tau+1,\phi;
\left( \begin{array}{c}
\vecs \\
\vecnull
\end{array}\right) 
+
\left( \begin{array}{cc}
1 & 1 \\
0 & 1
\end{array}\right)
\left( \begin{array}{c}\vecx \\ \vecy\end{array}\right) ,t 
+\tfrac{1}{2}\vecs\cdot\vecy\right) 
= \Theta_f\left(\tau,\phi;\left( \begin{array}{c}\vecx \\\vecy\end{array}\right),t\right),
$$
with
$$
\vecs=\trans{(\tfrac12,\tfrac12,\ldots,\tfrac12)}\in\RR^k.
$$
\endproclaim

\demo{Proof}
Clearly for any $f\in\Sw(\RR^k)$
$$
\sum_{\vecm\in\ZZ^k}
\left[\tilde R\left(\i+1,0;\left( \begin{array}{c}
\vecs \\
\vecnull
\end{array}\right) 
,0\right) f\right](\vecm)
=
\sum_{\vecm\in\ZZ^k}
f(\vecm),
$$
and hence also (replace $f$ with $\tilde R(\tau,\phi;\vecxi,t)f$)
$$
\sum_{\vecm\in\ZZ^k}
\left[\tilde R\left(\i+1,0;\left( \begin{array}{c}
\vecs \\
\vecnull
\end{array}\right) 
,0\right) \tilde R(\tau,\phi;\vecxi,t)f\right](\vecm)
=
\sum_{\vecm\in\ZZ^k}
[\tilde R(\tau,\phi;\vecxi,t) f](\vecm).
$$
We conclude by noticing
\begin{eqnarray*}
&& \hskip-.5in
\tilde R\left(\i+1,0;\left( \begin{array}{c}
\vecs \\
\vecnull
\end{array}\right) 
,0\right) \tilde R\left(\tau,\phi;\left( \begin{array}{c}\vecx \\\vecy\end{array}\right) ,t\right) 
\\
&&\qquad =
\tilde R\left(\tau+1,\phi;
\left( \begin{array}{c}
\vecs \\
\vecnull
\end{array}\right) 
+
\left( \begin{array}{cc}
1 & 1 \\
0 & 1
\end{array}\right)
\left( \begin{array}{c}\vecx \\\vecy\end{array}\right) ,t 
+\tfrac{1}{2}\vecs\cdot\vecy\right) ,
\end{eqnarray*} 
where we have used that $c((\i,0),(\tau,\phi))=1$
since $(\i,0)$ is an upper triangular matrix; cf.~3.6.
\enddemo

\nonumproclaim{Jacobi 3} 
$$
\Theta_f\left(\tau,\phi;\left( \begin{array}{c}\veck \\\vecl \end{array}\right)
+\vecxi,r+t
+\tfrac12\omega\left(\left( \begin{array}{c}\veck \\\vecl \end{array}\right),\vecxi\right)\right)
=(-1)^{\veck\cdot\vecl}\,
\Theta_f(\tau,\phi;\vecxi,t)
$$
for any $\veck,\vecl\in\ZZ^k${\rm ,} $r\in\ZZ$.
\endproclaim

\demo{Proof}
By virtue of~3.2 we have for all $f$
$$
\sum_{\vecm\in\ZZ^k}
\left[W\left(\left( \begin{array}{c}\veck \\\vecl \end{array}\right),r\right) f\right] (\vecm)
=
e(-\tfrac12 \veck\cdot\vecl)
\sum_{\vecm\in\ZZ^k}
f(\vecm), \pagebreak
$$
and therefore, replacing $f$ with $W(\vecxi,t)\tilde R(\tau,\phi)f$,
\begin{eqnarray*} &&\hskip-.5in
\sum_{\vecm\in\ZZ^k}
\left[W\left(\left( \begin{array}{c}\veck \\\vecl \end{array}\right),r\right) 
W(\vecxi,t)\tilde R(\tau,\phi)f\right] (\vecm)
\\ &&\qquad\qquad =
e(-\tfrac12 \veck\cdot\vecl)
\sum_{\vecm\in\ZZ^k}
[W(\vecxi,t)\tilde R(\tau,\phi)f](\vecm) ,
\end{eqnarray*}
which gives the desired result.
\enddemo

4.5.  In what follows, we shall only need to consider products of theta sums
of the form
$$
\Theta_f(\tau,\phi;\vecxi,t)\overline{\Theta_g(\tau,\phi;\vecxi,t)} ,
$$
where $f,g\in\Sw(\RR^k)$. Clearly such combinations do not depend
on the $t$-variable. Let us therefore define the semi-direct product
group
$$
G^k=\SL(2,\RR)\ltimes \RR^{2k}
$$
with multiplication law
$$
(M;\vecxi)(M';\vecxi') =(MM';\vecxi+M\vecxi') ,
$$
and put
$$
\Theta_f(\tau,\phi;\vecxi) =
v^{k/4} 
\sum_{\vecm\in\ZZ^k} f_\phi( (\vecm-\vecy) v^{1/2}) 
e(\tfrac12\|\vecm-\vecy\|^2 u + \vecm\cdot\vecx) .
$$
By virtue of Lemma 4.3 and the Iwasawa parametrization
3.8,
$\Theta_f\overline{\Theta_g}$ is a continuous $\CC$-valued function
on $G^k$.

\vglue6pt 4.6. A short calculation yields that the set
$$
\Gamma^k=\left\{ 
\left( \left( \begin{array}{cc}
a & b \\
c & d
\end{array}\right)
; 
\left( \begin{array}{c}
ab \vecs \\
cd \vecs
\end{array}\right)
+\vecm \right)
: \; \left( \begin{array}{cc}
a & b \\
c & d
\end{array}\right) \in\SL(2,\ZZ),\; \vecm\in\ZZ^{2k}
\right\} \subset G^k,
$$
with
$\vecs=\trans{(\tfrac12,\tfrac12,\ldots,\tfrac12)}\in\RR^k$,
is closed under multiplication and inversion, and therefore 
forms a subgroup of $G^k$.
Note also that the subgroup 
$$
N=\{1\}\ltimes\ZZ^{2k}
$$ 
is normal in $\Gamma^k$.

\nonumproclaim{4.7. Lemma} 
 $\Gamma^k$ is generated by the elements
$$
\left( \left( \begin{array}{cc}
0 & -1 \\
1 & 0
\end{array}\right)
; \vecnull \right), \quad
\left( \left( \begin{array}{cc}
1 & 1 \\
0 & 1
\end{array}\right)
; \left( \begin{array}{c}
\vecs \\
\vecnull
\end{array}\right)\right), \quad
\left( \left( \begin{array}{cc}
1 & 0 \\
0 & 1
\end{array}\right)
; \vecm\right), \quad \vecm\in\ZZ^{2k} . \pagebreak
$$
\endproclaim

\demo{Proof}
The map
$$
\SL(2,\ZZ)\rightarrow N \backslash\Gamma^k,
\qquad
\left( \begin{array}{cc}
a & b \\
c & d
\end{array}\right)
\mapsto
\left(\left( \begin{array}{cc}
a & b \\
c & d
\end{array}\right)
; 
\left( \begin{array}{c}
ab \vecs \\
cd \vecs
\end{array}\right)
+\ZZ^{2k}\right)
$$
defines a group isomorphism.
The matrices 
\def\nhs{\hskip-3pt}
$
( \begin{array}{cc}
\scr{0} \nhs&\nhs \scr{-1} \\[-4pt]
\scr{1} \nhs&\nhs \scr{0}
\end{array})
$  
and
$(\begin{array}{cc}
\scr{1} \nhs &\nhs \scr{1} \\[-4pt]
\scr{0} \nhs &\nhs  \scr{1}
\end{array})$
generate $\SL(2,\ZZ)$, hence the lemma.
\enddemo

\nonumproclaim{4.8. Proposition}
The left action of the group $\Gamma^k$ on $G^k$ 
is properly discontinuous.
A fundamental domain of $\Gamma^k$ in $G^k$ is given by
$$
{\cal F}_{\Gamma^k} = {\cal F}_{\SL(2,\ZZ)}
\times \{ \phi\in[0,\pi) \} \times \{\vecxi\in[-\tfrac12,\tfrac12)^{2k} \}.
$$
where ${\cal F}_{\SL(2,\ZZ)}$ 
is the fundamental domain in $\H$ of the modular group
$\SL(2,\ZZ)${\rm ,} given by $\{ \tau\in\H: u\in[-\tfrac12,\tfrac12), |\tau|> 1 \}$.
\endproclaim

\demo{Proof}
As mentioned before, the matrices 
$
( \begin{array}{cc}
\scr{0}\nhs &\nhs \scr{-1} \\[-4pt]
\scr{1}\nhs & \scr{0}
\end{array})
$
and
$(\begin{array}{cc}
\scr{1}\nhs & \nhs \scr{1} \\[-4pt]
\lscr{0} & \rscr{1}
\end{array})$
generate $\SL(2,\ZZ)$, which explains ${\cal F}_{\SL(2,\ZZ)}$. 
Note furthermore that 
$
( \begin{array}{cc}
\lscr{-1} & \rscr{0} \\[-4pt]
\lscr{0} & \rscr{-1}
\end{array})
$
generates the shift $\phi\mapsto\phi+\pi$. 
\enddemo

\nonumproclaim{4.9. Proposition}
 For $f,g\in\Sw(\RR^k)${\rm ,}
$\Theta_f(\tau,\phi;\vecxi)\overline{\Theta_g(\tau,\phi;\vecxi)}$ 
is invariant under the left action
of $\Gamma^k$.
\endproclaim

\demo{Proof}
This follows directly from Jacobi~1--3,
since the left action of the generators from 4.7 is
$$
\left(\tau,\phi;\left( \begin{array}{c}\vecx \\ \vecy\end{array}\right)\right)
\mapsto \left(-\frac{1}{\tau},\phi+\arg\tau;
\left( \begin{array}{c}-\vecy \\ \vecx\end{array}\right)\right) ,
$$
$$
(\tau,\phi;\vecxi)\mapsto\left(\tau+1,\phi;
\left( \begin{array}{c}
\vecs \\
\vecnull
\end{array}\right) 
+
\left( \begin{array}{cc}
1 & 1 \\
0 & 1
\end{array}\right)
\left( \begin{array}{c}\vecx \\\vecy\end{array}\right)\right),
$$
and
$$
(\tau,\phi;\vecxi)\mapsto
(\tau,\phi;\vecxi+\vecm) ,
$$
respectively.
\enddemo

We find the following uniform estimate.

\nonumproclaim{4.10 Proposition} Let $f,g\in\Sw(\RR^k)$.  For any $R>1${\rm ,} 
\begin{eqnarray*}
&&\hskip-.5in\Theta_f\left(\tau,\phi;\left( \begin{array}{c}\vecx \\ \vecy \end{array}\right)\right)
\overline{\Theta_g\left(\tau,\phi;\left( \begin{array}{c}\vecx \\ \vecy \end{array}\right)\right)}
\\
&&\qquad = v^{k/2} 
\sum_{\vecm\in\ZZ^k}
f_\phi((\vecm-\vecy) v^{1/2}) \overline{g_\phi((\vecm-\vecy) v^{1/2})}
+O_R(v^{-R})
\end{eqnarray*}
uniformly for all $(\tau,\phi;\vecxi)\in G^k$ with $v>\frac12$.
In addition{\rm ,}
\begin{eqnarray*}
&&\hskip-.5in 
\Theta_f\left(\tau,\phi;\left( \begin{array}{c}\vecx \\ \vecy \end{array}\right)\right)
\overline{\Theta_g\left(\tau,\phi;\left( \begin{array}{c}\vecx \\ \vecy \end{array}\right)\right)}
\\
&&\qquad = v^{k/2} 
f_\phi((\vecn-\vecy) v^{1/2}) \overline{g_\phi((\vecn-\vecy) v^{1/2})}
+O_R(v^{-R}),
\end{eqnarray*}
uniformly for all $(\tau,\phi;\vecxi)\in G^k$ with $v>\frac12${\rm ,}
$\vecy\in\vecn+[-\tfrac12,\tfrac12]^k$ and $\vecn\in\ZZ^k$.
\endproclaim

\demo{Proof}
Suppose $\vecy\in\vecn+[-\frac12,\frac12]^k$ for an arbitrary integer 
$\vecn\in\ZZ^k$.

By virtue of Lemma 4.3 we have for any $T>1$
$$
|f_\phi((\vecm-\vecy)v^{1/2})| \leq c_T (1+\|\vecm-\vecy\|v^{1/2})^{-T} 
=O_T(\|\vecm-\vecn\|^{-T}\, v^{-T/2}),
$$
which holds uniformly for $v>\frac12$, 
$\phi\in\RR$ and $\vecy\in\vecn+[-\frac12,\frac12]^k$,
if $\vecm\neq\vecn$. 

Likewise for $g_\phi$,
$$
 |g_\phi((\tilde{\vecm}-\vecy)v^{1/2})| \leq \tilde c_T 
(1+\|\tilde{\vecm}-\vecy\|v^{1/2})^{-T} 
=O_T(\|\tilde{\vecm}-\vecn\|^{-T}\, v^{-T/2}),
$$
again uniformly for $v>\frac12$, 
$\phi\in\RR$ and $\vecy\in\vecn+[-\frac12,\frac12]^k$,
if $\tilde{\vecm}\neq\vecn$. 

Hence the leading order contributions come from terms with 
$\tilde{\vecm}=\vecm$, the sum of all other terms
contributes $O_T(v^{-T/2})$.
\enddemo

The following lemmas will be useful later on.

\nonumproclaim{4.11 Lemma} 
 The subgroup
$$
\Gamma_\theta \ltimes \ZZ^{2k},
$$
where 
$$
\Gamma_\theta=\left\{ 
\left( \begin{array}{cc}
a & b \\
c & d
\end{array}\right) \in\SL(2,\ZZ):\; ab\equiv cd \equiv 0\bmod 2
\right\} 
$$ 
denotes the theta group{\rm ,}  is  of index three in $\Gamma^k$.
\endproclaim

\demo{Proof}
It is well known \cite{Hejhal83} that
$\Gamma_\theta$ is of index three in $\SL(2,\ZZ)$
and
$$
\SL(2,\ZZ) = \bigcup_{j=0}^2 \Gamma_\theta 
\left( \begin{array}{cc} 0 & -1 \\ 1 & 1 \end{array}\right)^j .
$$
By virtue of the group isomorphism employed in the proof of Lemma
4.7, we infer that
\vglue12pt
\hfill ${\displaystyle
\Gamma^k = \bigcup_{j=0}^2 (\Gamma_\theta\ltimes\ZZ^{2k}) 
\left(\left( \begin{array}{cc} 0 & -1 \\ 1 & 1 \end{array}\right); 
\left( \begin{array}{c} \vecnull \\ \vecs \end{array}\right)\right)^j .
}$
\enddemo
\vglue12pt

\nonumproclaim{4.12 Lemma}  $\Gamma^k$ is of finite index in $\SL(2,\ZZ)\ltimes (\frac12\ZZ)^{2k}$.
\endproclaim

\demo{Proof}
The subgroup $\Gamma_\theta\ltimes\ZZ^{2k}\subset\Gamma^k$ 
is of finite index in
$\SL(2,\ZZ)\ltimes \ZZ^{2k}$ and thus also in
$\SL(2,\ZZ)\ltimes (\frac12\ZZ)^{2k}$.
\enddemo

 4.13. {\it Remark}. 
Note that 
$$
\SL(2,\ZZ)\ltimes (\tfrac12\ZZ)^{2k}
=
\left(\left( \begin{array}{cc} 
\tfrac12 & 0 \\ 0 & \tfrac12
\end{array}\right); \vecnull\right)
(\SL(2,\ZZ)\ltimes \ZZ^{2k})
\left(\left( \begin{array}{cc} 
2 & 0 \\ 0 & 2
\end{array}\right); \vecnull\right),
$$
i.e., $\SL(2,\ZZ)\ltimes (\frac12\ZZ)^{2k}$ is isomorphic to
$\SL(2,\ZZ)\ltimes\ZZ^{2k}$.
\vglue8pt

4.14. In this paper, we will be interested in the case
of quadratic forms in two variables, i.e., $k=2$.
The corresponding theta sum (defined for general $k$ in 4.5)
reads then
\begin{eqnarray*}
\Theta_f(\tau,\phi;\vecxi) &= &v^{1/2} 
\sum_{(m,n)\in\ZZ^2} 
f_\phi( (m-y_1) v^{1/2},(n-y_2) v^{1/2}) \\
&& \times\ e(\tfrac12(m-y_1)^2 u + \tfrac12(n-y_2)^2 u + mx_1 +nx_2) ,
\end{eqnarray*}
where $\vecxi=\trans(x_1,x_2,y_1,y_2)\in\RR^4$.
This theta sum is related to the one introduced in
Section \ref{secplan} by
\begin{eqnarray*}
&&\qquad \theta_{\psi_1}(u,\lambda)\overline{\theta_{\psi_2}(u,\lambda)}
=\Theta_f(\tau,\phi;\vecxi)\overline{\Theta_g(\tau,\phi;\vecxi)}
\\
\noalign{\noindent 
with}
&&
\quad\tau=u+\i\frac1\lambda, \qquad \phi=0, \qquad
\vecxi=\trans(0,0,\alpha,\beta),\\
\noalign{\noindent 
and}
&&
f(w_1,w_2)=\psi_1(w_1^2+w_2^2),\qquad g(w_1,w_2)=\psi_2(w_1^2+w_2^2).
\end{eqnarray*}
Recall that $f_\phi|_{\phi=0}=f$ and likewise $g_\phi|_{\phi=0}=g$.

The crucial advantage in dealing with $\Theta_f$ rather than the
original $\theta_\psi$ is that the extra set of variables allows
us to realize $\Theta_f$ as a function on a finite-volume
manifold and to employ ergodic-theoretic techniques.

 \vglue-8pt
\section{Unipotent flows}\label{secunipotent}
\vglue-12pt
5.1.
Put
$$
\Psi_0^t=
\left(\left( \begin{array}{cc}
1 & t \\
0 & 1
\end{array}\right) ; \vecnull \right).
$$
For $t\in\RR$, 
$\Psi_0^t$ generates a unipotent one-parameter-subgroup of $G^k$,
denoted by $\Psi_0^\RR$.
For any lattice $\Gamma$ in $G^k$,
we now define the flow $\Psi^t: 
\Gamma\backslash G^k\rightarrow\Gamma\backslash G^k$
by right translation by $\Psi_0^t$,
$$
\Psi^t(g) := g \Psi_0^t .
$$
Hence for $g=(M;\vecxi)$ we have
$$
\Psi^t(g)=\left(M \left( \begin{array}{cc}
1 & t \\
0 & 1
\end{array}\right);\vecxi\right).
$$ 
When projected onto $\Gamma\backslash\SL(2,\RR)$,
this flow becomes the classical  horocycle flow. \pagebreak

  5.2. Similarly,  
$$
\Phi_0^t=
\left(\left( \begin{array}{cc}
\e^{-t/2} & 0 \\
0 & \e^{t/2}
\end{array}\right) ; \vecnull \right),
$$
generates a one-parameter-subgroup of $G^k$.
The flow $\Phi^t: \Gamma\backslash G^k\rightarrow\Gamma\backslash G^k$
defined by
$$
\Phi^t(g) := g \Phi_0^t ,
$$
represents a lift of the classical geodesic flow on 
$\Gamma\backslash\SL(2,\RR)$.

\vglue12pt 5.3.
We are interested in averages of the form
$$
\int F(u+\i v,0;\vecxi) \; h(u)\, du
$$
where $F$ is a continuous function $\Gamma\backslash G^k\rightarrow \RR$,
and $h$ is a continuous probability density 
with compact support.
Setting $g_0=(\i,0;\vecxi)$, and $v=\e^{-t}$, we may write the above
integral as
$$
\rho_t(F)=\int F(g_0 \Psi_0^u \Phi_0^t) \; h(u)\, du =
\int F\circ\Phi^t\circ\Psi^u(g_0) \; h(u)\, du ,
$$
which may therefore be interpreted as the average along an orbit
of the unipotent flow $\Psi^u$, which is translated by $\Phi^t$.
Since $\rho_t(1)=1$, 
$\rho_t$ defines a probability measure on $\Gamma\backslash G^k$.

\nonumproclaim{5.4. Proposition} Let $\Gamma$ be a subgroup of $\SL(2,\ZZ)\ltimes\ZZ^{2k}$ of finite
index. Then the family of probability measures $\{\rho_t: \; t\geq 0\}$
is relatively compact{\rm ,} i.e.{\rm ,} every sequence of measures contains a
subsequence which converges weakly to a probability measure
on $\Gamma\backslash G^k$.
\endproclaim

\demo{Proof}
Consider the function
$$
X_R(\tau)= 
\sum_{\gamma\in\{\Gamma_\infty\cup(-1)\Gamma_\infty\}\backslash\SL(2,\ZZ)}
\chi_R(\Im(\gamma\tau)),
$$
where $\chi_R$ is the characteristic function of the open interval
$(R,\infty)$, and
$$
\Gamma_\infty = \left\{ \left( \begin{array}{cc} 1 & m \\ 0 & 1 \end{array}\right):
m\in\ZZ\right\} \subset \SL(2,\ZZ) .
$$
For $u+\i v\in{\cal F}_{\SL(2,\ZZ)}$, we thus have
$$
X_R (u+\i v) =
\left\{ \begin{array}{ll} 
1 & (v> R) \\
0 & (v\leq R).
\end{array}\right. 
$$
Because $\Gamma$ is a finite index subgroup of $\SL(2,\ZZ)\ltimes\ZZ^{2k}$,
$X_R$ represents the characteristic function of a set in 
$\Gamma\backslash G^k$, whose complement is compact.

By construction,
the function $X_R$ is independent of $\phi$ and $\vecxi$;
we can therefore apply the
equidistribution theorem for arcs of long closed horocycles
on $\Gamma\backslash\H$
(see, e.g., \cite{Hejhal95} and \cite[Cor.~5.2]{Marklof98}),
which yields for $g_0=(\i,0;\vecxi)$,
$$
\lim_{t\rightarrow\infty}\rho_t(X_R)
=\lim_{v\rightarrow 0} 
\int X_R(u+\i v) \; h(u)\, du 
= \frac{1}{\mu({{\cal F}_{\SL(2,\ZZ)}})}
\int_{{\cal F}_{\SL(2,\ZZ)}} X_R(u+\i v) \, \frac{du\,dv}{v^2}.
$$
Now
$$
\int_{{\cal F}_{\SL(2,\ZZ)}} X_R(u+\i v) \, \frac{du\,dv}{v^2} 
= \int_R^\infty \frac{dv}{v^2}
= R^{-1} .
$$
Hence, given any $\varepsilon>0$, we find some $R>1$ such that
$$
\sup_{t\geq 0} \rho_t(X_R) \leq \varepsilon. 
$$
The family of $\rho_t$ is therefore tight, 
and the proposition follows from the Helly-Prokhorov theorem
\cite{Shiryaev95}.
\enddemo

\nonumproclaim{5.5. Proposition} If $\nu$ is a weak limit of a subsequence of the
probability measures $\rho_t$ with $t\rightarrow\infty${\rm ,} then $\nu$ is
invariant under the action of $\Psi^\RR${\rm ,} i.e.{\rm ,} $\nu\circ\Psi^\RR=\nu$.
\endproclaim

\demo{Proof}
Suppose $\{\rho_{t_i}:\,i\in\NN\}$ is a convergent subsequence
with weak limit~$\nu$. That is, for any bounded continuous function
$F$, we have
$$
\lim_{i\rightarrow\infty}\rho_{t_i}(F)=\nu(F).
$$
For any fixed $s\in\RR$, we find
\begin{eqnarray*}
\rho_{t}(F\circ\Psi^s)
&=&\int F(g_0 \Psi_0^u \Phi_0^{t}\Psi_0^s) \; h(u)\, du 
=\int F(g_0 \Psi_0^{u+s\exp(-t)} \Phi_0^{t}) \; h(u)\, du \\
&=&\int F(g_0 \Psi_0^{u} \Phi_0^{t}) \; h(u-s\exp(-t))\, du .
\end{eqnarray*}
Furthermore
\begin{eqnarray*}
\big|\rho_{t}(F\circ\Psi^s)-\rho_{t}(F)\big| 
&=&\big| \int F(g_0 \Psi_0^{u} \Phi_0^{t}) \; \big[h(u-s\exp(-t))
-h(u)\big]\, du \big|   \\
&\leq& (\sup |F|)  
\int \big|h(u-s\exp(-t))
-h(u)\big|\, du .
\end{eqnarray*}
Hence, given any $\varepsilon>0$, we find a $T$ such that 
$$
|\rho_{t}(F\circ\Psi^s)-\rho_{t}(F)|<\varepsilon
$$
for all $t>T$.
Because the function $\tilde F=F\circ\Psi^s$ ($s$ is fixed) is bounded
continuous, the limit
$$
\lim_{i\rightarrow\infty}\rho_{t_i}(F\circ\Psi^s)=\nu(F\circ\Psi^s)\pagebreak
$$
exists, and we know from the above inequality that
$$
|\nu(F\circ\Psi^s)-\nu(F)|\leq\varepsilon
$$
for any $\varepsilon>0$. Therefore $\nu(F\circ\Psi^s)=\nu(F)$.
\enddemo

5.6. Ratner \cite{Ratner91}, \cite{Ratner91b} gives a classification of
all ergodic $\Psi^\RR$-invariant measures on $\Gamma\backslash G^k$.
We will now investigate which of these measures are possible
limits of the sequence $\{\rho_t\}$. The answer will be unique,
translates of orbits of $\Psi^\RR$ become equidistributed. 

\nonumproclaim{5.7. Theorem} Let $\Gamma$ be a subgroup of $\SL(2,\ZZ)\ltimes\ZZ^{2k}$ of
finite index. Fix some point 
$$
g_0=\left(\i,0;\left( \begin{array}{c}\vecnull \\ \vecy \end{array}\right)\right)
\in\Gamma\backslash G^k
$$ 
such that the components of the vector 
$(\trans\vecy,1)\in\RR^{k+1}$ are linearly
independent over $\QQ$.
Let $h$ be a continuous probability density $\RR\rightarrow\RR_+$
with compact support.
Then{\rm ,} for any bounded continuous function $F$ on 
$\Gamma\backslash G^k${\rm ,} 
$$
\lim_{t\rightarrow\infty}
\int_\RR F\circ\Phi^t\circ\Psi^u(g_0) \; h(u)\, du
= \frac{1}{\mu(\Gamma\backslash G^k)}\int_{\Gamma\backslash G^k} F \, d\mu
$$
where $\mu$ is the Haar measure of $G^{k}$.
\endproclaim

This theorem is a special case of Shah's more 
general Theorem 1.4 in \cite{Shah96}
on the equidistribution of translates of unipotent orbits. 
Because of the simple structure
of the Lie groups studied here, the proof of Theorem~5.7
is less involved than in the general context.

\vglue8pt
5.8. Before we begin with the proof of Theorem~{\rm 5.7,}
we consider the special test function 
$$
F_\delta(M; \vecxi)
=\sum_{\gamma\in\SL(2,\ZZ)} 
f_\delta(\gamma M)\; \eta_D(\gamma\vecxi),
$$
with {\rm (}\/in the Iwasawa parametrization {\rm 3.8)}
$$
f_\delta(M)=f_\delta(\tau,\phi)=\chi_1(u +v \cot\phi)\;
\chi_2(v^{-1/2}\cos\phi)\;
\chi_3(v^{-1/2}\sin\phi)
$$
where $\chi_j$ ($j=1,2,3$) is the characteristic function
of the interval $[s_j,s_j+\delta_j]$. We assume
in the following that $s_j$ ranges over the fixed compact interval $I_j${\rm ,}
and that $I_3$ is furthermore properly contained in $\RR^+$, i.e.,
$s_3\ge\overline s$ for some constant $\overline s>0$.
Clearly $f_\delta$ has compact support in $\SL(2,\RR)$.
The function $\eta_D:\TT^{2k}\rightarrow\RR$ is the characteristic function
of a domain $D$ in $\TT^{2k}$ with smooth boundary.

Clearly, $F_\delta$ may be viewed as a function on 
$\Gamma\backslash G^k$, for $\Gamma$ is a subgroup
of $\SL(2,\ZZ)\ltimes\ZZ^{2k}$.

\nonumproclaim{5.9. Lemma}  Suppose the components of the vector 
$(\trans\vecy,1)\in\RR^{k+1}$ are linearly
independent over $\QQ$. Then{\rm ,} given intervals $I_1,I_2,I_3$ as above{\rm ,}
there exists a constant $C>0$ such that{\rm ,} for 
any domain $D\subset\TT^{2k}$ with smooth boundary{\rm ,} 
$\delta_1,\delta_2,\delta_3>0$ {\rm (}\/sufficiently small\/{\rm )}
and $s_1\in I_1,s_2\in I_2,s_3\in I_3${\rm ,}  
$$
\limsup_{v\rightarrow 0}\int_\RR F_\delta\left(u+\i v, 0;
\left( \begin{array}{c}\vecnull \\ \vecy \end{array}\right)\right) \; h(u)\, du
\leq C \delta_1\delta_2(s_3+\delta_3) \int_{\TT^{2k}} \eta_D(\vecxi)d\vecxi .
$$ 
The constant $C$ may depend on the choice of $h,\vecy,I_1,I_2,I_3$.
\endproclaim

5.10. {\it Proof.}
\vglue4pt
5.10.1. Given any $\varepsilon>0$ and
any domain $D\subset\TT^{2k}$ with smooth boundary, 
we can cover $D$ by a large but finite number of 
nonoverlapping cubes $C_j\subset\TT^{2k}$, in such
a way that 
$$
\eta_D \leq \sum_j \eta_{C_j}, \qquad 
\int_{\TT^{2k}} \left(\sum_j\eta_{C_j}-\eta_D\right) d\vecxi <\varepsilon .
$$
We may therefore assume without loss of generality that 
$\eta_D(\vecxi)$ is the characteristic function of an arbitrary cube
in $\TT^{2k}$, i.e.,
$\eta_D(\vecxi)=\eta_1(\vecx)\eta_2(\vecy)$,
where $\eta_1,\eta_2$ are characteristic functions of 
arbitrary cubes in $\TT^k$.

\vglue8pt 5.10.2.
We recall that for
$\gamma=(\begin{array}{cc} \lscr{a} & \rscr{b} \\[-4pt] \lscr{c} & \rscr{d} \end{array})$,
$$
F_\delta\left(u+\i v,0;\left( \begin{array}{c}\vecnull \\ \vecy \end{array}\right)\right)
=\sum_\gamma f_\delta\left(\frac{a(u+\i v)+b}{c(u+\i v)+d},\arg(c\tau+d)\right) \;
\eta_1(b\vecy)\,\eta_2(d \vecy) .
$$
In particular (with $\phi=0$),
$$
v_\gamma^{-1/2}\cos\phi_\gamma=v^{-1/2}(cu+d), \qquad
v_\gamma^{-1/2}\sin\phi_\gamma=v^{1/2} c,
$$  
$$
u_\gamma=\Re\frac{a(u+\i v)+b}{c(u+\i v)+d}
=\frac{a}{c}-\frac{1}{c}\,\frac{cu+d}{|c\tau+d|^2}
=\frac{a}{c}-v_\gamma\cot\phi_\gamma.
$$
One then finds that
$$
F_\delta\left(u+\i v,0;\left( \begin{array}{c}\vecnull \\ \vecy \end{array}\right)\right)
=\sum_\gamma 
\chi_1(\frac{a}{c})\;
\chi_2(v^{-1/2}(cu+d))\;
\chi_3(cv^{1/2})\; 
\eta_1(b\vecy)\;
\eta_2(d\vecy) ,
$$
which, after being integrated against $h(u)du$, yields
$$
v \sum_\gamma 
\chi_1(\frac{a}{c})\;
\chi_3(cv^{1/2})\; 
\eta_1(b\vecy)\;
\eta_2(d\vecy)
\int\chi_2(cv^{1/2}t)\;h(vt-\frac{d}{c}) \; dt.
$$
The compactness of the support of $h$ implies that 
$\frac{d}{c}=vt+O(1)$, and \pagebreak hence $|d|\ll |s_2+\delta_2| v^{1/2} + |c|$,
i.e., $|d| \ll |c|$ for $v$ small.
Therefore
\begin{eqnarray*}
&&\hskip-.5in\int F_\delta\left(u+\i v,0;\left( \begin{array}{c}\vecnull \\ \vecy \end{array}\right)\right)\, h(u)du
\\
&&\ll 
 \delta_2\, v^{1/2} \sum_{{\gamma\atop |d|\leq A|c|}} \frac{1}{|c|}
\chi_1\left(\frac{a}{c}\right)\;
\chi_3(cv^{1/2})\; 
\eta_1(b\vecy)\;
\eta_2(d\vecy) ,
\end{eqnarray*}
where $A>0$ 
and the implied constant depend only on $h$, if $v$ is small enough.

\vglue8pt 5.10.3. 
There are only finitely many terms with $d=0$, which thus give a total 
contribution of order $v^{1/2}$; we will thus assume in the following
$d\neq 0$. Likewise, if $b=0$, we have $ad=1$ and $c\in\ZZ$. This leads
to a contribution of order $v^{1/2} |\log v|$, which tends to zero
in the limit $v\rightarrow 0$.

The solutions of the equation $ad-bc=1$ with $b,d\neq 0$ can be obtained
in the following way. Take nonzero coprime integers $b,d\in\ZZ$, 
$\gcd(b,d)=1$,
and suppose $a_0,c_0$ solves $a_0d-bc_0=1$. (Such a solution can always be
found.) All other solutions must then be of the form $a=a_0+mb$,
$c=c_0+md$ with $m\in\ZZ$. We may assume without loss of generality
that $0\leq c_0\leq |d|-1$. So, for $v$ sufficiently small,
\begin{eqnarray*}
 &&\hskip-.75in \int F_\delta\left(u+\i v,0;\left( \begin{array}{c}\vecnull \\ \vecy \end{array}\right)\right)\, 
h(u)du\\
&\ll &
 \delta_2\, v^{1/2} \sum_{{b,d,m\in\ZZ\atop 0<|d|\leq A|c_0+md|}} 
\frac{1}{|c_0+md|}
\chi_1\left(\frac{b}{d}+\frac{1}{(c_0+md)d}\right)\\
&&\times
\chi_3((c_0+md)v^{1/2})\; 
\eta_1(b\vecy)\;
\eta_2(d\vecy) +O_{\delta,\eta}(v^{1/2}\log v),
\end{eqnarray*}
where $a_0=a_0(b,d)$ and $c_0=c_0(b,d)$ are chosen as above.
We have dropped the restriction that $\gcd(b,d)=1$.

For terms with $|m|>1$, we obtain upper bounds by observing
$$
\frac{1}{|c_0+md|}\leq \frac{1}{(|m|-1)|d|},
$$
and replacing the restriction imposed by $\chi_3$ with the condition 
$(|m|-1)|d|\leq v^{-1/2}(s_3+\delta_3)$. For terms with $m=0,\pm 1$,
we have 
$$
\frac{1}{|c_0+md|} \leq \frac{A}{|d|}
$$
and we replace the restriction corresponding to $\chi_3$ with
$|d|\leq A v^{-1/2}(s_3+\delta_3)$, since $|d|\leq A|c_0+md|$.

The restriction coming from $\chi_1$ means for $d>0$ 
$$
s_1 d -\frac{1}{c_0+md} \leq b \leq (s_1+\delta_1)d-\frac{1}{c_0+md},
$$
which we extend to
$$
s_1 d -\frac{A}{|d|} \leq b \leq (s_1+\delta_1)d +\frac{A}{|d|},
$$
and for $d<0$,
$$
(s_1+\delta_1) d -\frac{1}{c_0+md} \leq b \leq  s_1 d-\frac{1}{c_0+md} ,
$$
which we extend to
$$
(s_1+\delta_1) d -\frac{A}{|d|} \leq b \leq  s_1 d + \frac{A}{|d|}.
$$

We thus have (with $n=|m|-1$ for $|m|>1$, and $n=1$ for $m=0,\pm 1$)
\begin{eqnarray*}
&&\hskip-.5in
\int F_\delta\left(u+\i v,0;\left( \begin{array}{c}\vecnull \\ \vecy \end{array}\right)\right)\, 
h(u)du\\
&&\ll 
 \delta_2\, v^{1/2} \sum_{b,d,n\in\ZZ} 
\frac{1}{|nd|}
\eta_1(b\vecy)\;
\eta_2(d\vecy) +O_{\delta,\eta}(v^{1/2}\log v),
\end{eqnarray*}
with the summation restricted to
$$
s_1 |d| -\frac{A}{|d|} \leq \pm b 
\leq (s_1+\delta_1)|d| + \frac{A}{|d|}, \enspace
n |d|\leq \max(A,1) v^{-1/2}(s_3+\delta_3), \enspace
n>0 .
$$

\vglue8pt 5.10.4.
Since the components of $(\trans\vecy, 1)$ are linearly independent
over $\QQ$,
Weyl's equidistribution theorem (\cite[Satz 4]{Weyl16})
implies that
$$
\sum_{s_1 |d| -\frac{A}{|d|} \leq \pm b 
\leq (s_1+\delta_1)|d| + \frac{A}{|d|}} 
\eta_1(b\vecy) \ll |d| \delta_1 \int_{\TT^k}\eta_1(\vecx)d\vecx,
$$
uniformly for $|d|>v^{-1/4}$ large enough.
For $|d|\leq v^{-1/4}$ we use the trivial bound
$$
\sum_{s_1 |d| -\frac{A}{|d|} \leq \pm b 
\leq (s_1+\delta_1)|d| + \frac{A}{|d|}} 
\eta_1(b\vecy) = O_{\delta,\eta}( v^{-1/4}) ,
$$
for small enough $v$.
Therefore
\begin{eqnarray*}
&&\hskip-.25in\int F_\delta\left(u+\i v,0;\left( \begin{array}{c}\vecnull \\ \vecy \end{array}\right)\right)\, 
h(u)du \\
&&  \ll  
\delta_1 \delta_2\, v^{1/2} \sum_{{n>0\atop
|d|\ll n^{-1} v^{-1/2}(s_3+\delta_3)} }
\frac{1}{n}
\eta_2(d\vecy)\; \int_{\TT^k} \eta_1(\vecx) d\vecx + 
O_{\delta,\eta}(v^{1/4} (\log v)^2) ,
\end{eqnarray*}
where the last term includes all contributions from terms with 
$|d|\leq v^{-1/4}$.

\vglue8pt 5.10.5.
We split the remaining sum over $n$ into terms with $0<n<v^{-1/4}$ and terms
with $n\geq v^{-1/4}$.
In the first case we have, for $v\to 0$,
$$
n v^{1/2} \sum_{0<|d|\ll n^{-1} v^{-1/2}(s_3+\delta_3)}
\eta_2(d\vecy) \ll (s_3+\delta_3)
\int_{\TT^k} \eta_2(\vecx) d\vecx 
$$
by Weyl's equidistribution theorem. For $n\geq v^{-1/4}$ one simply 
uses the trivial bound
$$
\sum_{0<|d|\ll n^{-1} v^{-1/2}(s_3+\delta_3)}
\eta_2(d\vecy) 
\ll  n^{-1}v^{-1/2} (s_3+\delta_3) .
$$

\vglue8pt 5.10.6.
We conclude
\begin{eqnarray*}
&&\hskip-18pt\limsup_{v\rightarrow 0}
\int F_\delta\left(u+\i v,0;\left( \begin{array}{c}\vecnull \\ \vecy \end{array}\right)\right)\, h(u)du
\\
&&\hskip-12pt  \ll
\delta_1 \delta_2 (s_3+\delta_3) \int_{\TT^k} \eta_1(\vecx) d\vecx  
\limsup_{v\rightarrow 0}
\left[ \sum_{n<v^{-1/4}} n^{-2} 
\int_{\TT^k} \eta_2(\vecx) d\vecx 
+ \sum_{n\geq v^{-1/4}} n^{-2} \right] .
\end{eqnarray*}
Since $\lim_{v\rightarrow 0} \sum_{n<v^{-1/4}} n^{-2}= \frac{\pi^2}{6}<\infty$ 
and $\lim_{v\rightarrow 0} \sum_{n\geq v^{-1/4}} n^{-2}=0$,
the lemma is proved.
\hfill\qed

\vglue12pt  5.11. {\it Proof of Theorem}~5.7.
\vglue8pt
5.11.1. 
By Propositions~5.4 and 5.5, we find a convergent subsequence
of $\rho_{t_i}$ with weak limit $\nu$ invariant under $\Psi^\RR$. 
Hence for any bounded continuous function $F$ on $\Gamma\backslash G^k$,
$$
\lim_{i\rightarrow\infty} \rho_{t_i}(F) = \nu(F).
$$

\vglue8pt
5.11.2.
Following \cite{Mozes95},
we denote by $\cal H$ the collection of all closed connected subgroups
$H$ of $G^k$ such that $\Gamma\cap H$ is a lattice in $H$ and
the subgroup, which is generated by all unipotent one-parameter subgroups
of $G^k$ contained in $H$, acts ergodically from the right
on $\Gamma\backslash\Gamma H$
with respect to the $H$-invariant probability measure. 
This collection is countable (\cite[Th.~1.1]{Ratner91}), and we call
${\cal H}^*\subset{\cal H}$ the set containing one representative
of each $\Gamma$-conjugacy class.

Because $\SL(2,\RR)\ltimes\{\vecnull\}$ and 
$\{1\}\ltimes \RR^{2k}$ are each generated by unipotent
one-parameter subgroups, so is $G^k$, which of course acts ergodically
(with respect to Haar measure $\mu$)
from the right on $\Gamma\backslash G^k$, and so $G^k\in {\cal H}$.

Let 
\begin{eqnarray*}
N(H)&=&\{g\in G^k: \Psi_0^\RR \subset g^{-1} H g\} ,
\\
S(H)&=&\bigcup_{H'\in{\cal H},\; H'\subset H,\; H'\neq H}
N(H') ,
\\
\noalign{\noindent and} 
T_H&=&\pi(N(H)\backslash S(H)),
\end{eqnarray*} 
where $\pi$ is the natural quotient map
$G^k\rightarrow \Gamma\backslash G^k$.
We denote by $\nu_H$ the restriction of $\nu$ on $T_H$.
Then, for any $g\in N(H)\backslash S(H)$, the group $g^{-1} H g$
is the smallest closed subgroup of $G^k$ which contains $\Psi_0^\RR$
and whose orbit through $\pi(g)$ is closed in $\Gamma\backslash G^k$
(cf.~\cite[Lemma 2.4]{Mozes95}).

For all Borel measurable subsets ${\cal A}\subset\Gamma\backslash G^k$,
the $\Psi^\RR$-invariant measure $\nu$
admits the decomposition (see \cite[Th.~2.2]{Mozes95})
$$
\nu({\cal A}) = \sum_{H\in{\cal H}^*} \nu_H({\cal A}) .
$$
Furthermore (see \cite{Mozes95} for details),
for any $\Psi^\RR$-ergodic component $\iota$ of $\nu_H$,
with $\iota$ a probability measure, there exists a $g\in N(H)$ such that 
$\iota$ is the unique\break $g^{-1}Hg$-right-invariant probability measure
on the closed orbit $\Gamma\backslash\Gamma H g$. 
In particular, if $\nu(\pi(S(G^k)))=0$, then $\nu=\mu$ (up to normalization).

\vglue8pt 5.11.3.
Let us suppose first that there is at least one 
$H\in{\cal H}$ with $\nu_H\neq 0$, 
whose projection onto the $\SL(2,\RR)$-component 
is a closed connected subgroup
$L$ of $\SL(2,\RR)$ with $L\neq\SL(2,\RR)$ (compare Appendix~B).
Let $\Lambda$ be the projection of 
$\Gamma$ onto its $\SL(2,\RR)$-component.
Since $\Gamma\cap H$ is a lattice in $H$, 
$\Lambda\cap L$ is a lattice in $L$.
We can therefore construct a bounded continuous function 
$F(\tau,\phi;\vecxi)=F(\tau,\phi)$ such that
$$
\int F d\nu  \neq \frac{1}{\mu(\Gamma\backslash G^k)}\, \int F d\mu .
$$
With $F$ independent of $\vecxi$,
we apply the equidistribution theorem for
long arcs of closed horocycles \cite{Hejhal95},  \cite{Marklof98}, which yields
$$
\lim_{t\rightarrow\infty} \rho_{t}(F) =
\frac{1}{\mu(\Gamma\backslash G^k)} \, \int F d\mu .
$$
For the above subsequence (5.11.1) we find, however,
$$
\lim_{i\rightarrow\infty} \rho_{t_i}(F) = \int F d\nu,
$$
which leads to a contradiction. We shall therefore assume in the
following that $L=\SL(2,\RR)$.

\vglue8pt 5.11.4. 
The most general form of a closed connected subgroup $H$,
for which $L=\SL(2,\RR)$
and which contains a conjugate of $\Psi_0^\RR$,
is (see Appendix~B)
$$
H= (1;\vecxi_0)H_0(1;-\vecxi_0), \quad H_0=\SL(2,\RR)\ltimes\Omega,
$$
where $\Omega$ is a closed connected subgroup of $\RR^{2k}$
(i.e., $\Omega$ is a closed linear subspace of $\RR^{2k}$),
which is invariant under the action of $\SL(2,\RR)$.
Since $\SL(2,\RR)\ltimes\{\vecnull\}$ and $\{1\}\ltimes\Omega$ are 
generated by unipotent one-parameter subgroups, the same holds for 
$H_0$ and hence for $H$. The right action of $H$
on $\Gamma\backslash\Gamma H$ is obviously ergodic with respect
to the (unique) $H$-invariant probability measure $\iota$,
and therefore  \pagebreak $H\in{\cal H}$.

 5.11.5.
Let us consider the orbit
$$
\Gamma\backslash\Gamma H g = 
\Gamma\backslash\Gamma (1;\vecxi_0)H_0 \tilde g 
$$
with $g\in N(H)$ and thus $\tilde g= (1;-\vecxi_0)g\in N(H_0)$.
Note that
$$
\left(1;\left( \begin{array}{c} \veca \\ \vecnull \end{array}\right)\right)\; \Psi_0^t\;
\left(1;\left( \begin{array}{c} \veca \\ \vecnull \end{array}\right)\right)^{-1}
=\Psi_0^t \in H_0
$$
for all $t\in\RR$, $\veca\in\RR^k$, and
$$
\left(1;\left( \begin{array}{c} \vecnull \\ \vecb \end{array}\right)\right) \;\Psi_0^t\;
\left(1;\left( \begin{array}{c} \vecnull \\ \vecb \end{array}\right)\right)^{-1}
=\left(1;\left( \begin{array}{c} -t \vecb \\ \vecnull \end{array}\right)\right) \Psi_0^t .
$$
The right-hand side is an element of $H_0$ for all $t\in\RR$ if and only
if 
$$
\left(\left( \begin{array}{cc} 0 & 1 \\ -1 & 0 \end{array}\right); \vecnull \right)
\left(1;\left( \begin{array}{c} -t \vecb \\ \vecnull \end{array}\right)\right)
\left(\left( \begin{array}{cc} 0 & -1 \\ 1 & 0 \end{array}\right); \vecnull \right)
=\left(1;\left( \begin{array}{c} \vecnull \\ t \vecb  \end{array}\right)\right) 
\in H_0 .
$$
We therefore have the explicit
representation
$$
N(H_0)=H_0 \; \left\{\left(1;\left( \begin{array}{c} \veca \\ \vecnull \end{array}\right)\right):
\;\veca\in\RR^k\right\} .
$$

\vglue8pt 5.11.6.
Let us suppose in the following that $\nu_H\neq 0$ for some $H\neq G^k$, i.e., 
$\Omega\neq \RR^{2k}$. 
We denote by $B_k(r)$ the open ball $\{\vecx\in\RR^k: \|\vecx\|<r\}$.
Then, for any $r>0$, we define
\begin{eqnarray*}
\Sigma(r)&=&\Gamma(1;\vecxi_0)H_0 
\left\{\left(1;\left( \begin{array}{c} \veca \\ \vecnull \end{array}\right)\right):
\;\veca\in B_k(r) \right\} \\
&=&
\left\{\left(M;\vecxi +M \left( \begin{array}{c} \veca \\ \vecnull \end{array}\right)\right):
\; M\in\SL(2,\RR),\; \vecxi\in\widetilde\Omega,\; \veca\in B_k(r) \right\} ,
\end{eqnarray*}
where $\widetilde\Omega= \Gamma(\vecxi_0+\Omega)$ is a closed subset
in $\RR^{2k}$.
We fix $r$ large enough so that the restriction of $\nu_H$
on $\Gamma\backslash\Sigma(r)$ is nonzero.

\vglue8pt 5.11.7. 
Let us discuss the structure of $\widetilde\Omega$ in more detail:
Since $\Gamma$ is of finite index in $\Gamma'=\SL(2,\ZZ)\ltimes\ZZ^{2k}$
we see that $\Gamma\cap H$
is of finite index in $\Gamma'\cap H$. Furthermore $\Gamma\cap H$ is a 
lattice in $H$, and so $\Gamma'\cap H$ is a lattice in $H$.
Then clearly $(1;-\vecxi_0)\Gamma'(1;\vecxi_0)\cap H_0$ 
must be a lattice in $H_0$.
With
$$
(1;-\vecxi_0)\Gamma'(1;\vecxi_0)=\{ (M;(M-1)\vecxi_0+\vecm):
\; M\in\SL(2,\ZZ),\;
\vecm\in\ZZ^{2k}\},
$$
the lattice property 
in turn implies that $(M-1)\vecxi_0\in \Omega+\ZZ^{2k}$ for all
$M$ in a finite index subgroup $\Lambda\subset\SL(2,\ZZ)$.
The orbit
$\SL(2,\ZZ)\vecxi_0/(\Omega+\ZZ^{2k})$ is therefore finite
 in $\RR^{2k}/(\Omega+\ZZ^{2k})$;
we denote by $\{\vecxi_0^{(1)},\vecxi_0^{(2)},\ldots,\vecxi_0^{(J)}\}$
a finite set of representatives.
With this, we conclude
$$
\Gamma'(\vecxi_0+\Omega)=\bigcup_{j=1}^J \vecxi_0^{(j)}+\Omega+\ZZ^{2k}.
$$
The fact that $(1;-\vecxi_0)\Gamma'(1;\vecxi_0)\cap H_0$ 
is a lattice in $H_0$ implies also that $\ZZ^{2k}\cap\Omega$ is a 
euclidean lattice
in $\Omega$. Hence there is a compact fundamental
domain ${\cal F}_{\ZZ^{2k}\cap\Omega}\subset\Omega$. We may therefore write
$$
\Gamma'(\vecxi_0+\Omega)=\bigcup_{j=1}^J \vecxi_0^{(j)}
+{\cal F}_{\ZZ^{2k}\cap\Omega}+\ZZ^{2k}.
$$
Note that ${\cal F}_{\ZZ^{2k}\cap\Omega}$
is also compact in $\RR^{2k}$, since $\Omega$ is closed. 

We conclude by observing that $\Gamma'(\vecxi_0+\Omega)$ is, 
of course, a finite covering of $\widetilde\Omega$, because $\Gamma$
has finite index in $\Gamma'$.

\vglue8pt 5.11.8.
Consider the subset $\Sigma_\delta(r)$ of $\Sigma(r)$,
given by
$$
\Sigma_\delta(r)
=
\Gamma 
\left\{\left(M;\vecxi +M \left( \begin{array}{c} \veca \\ \vecnull \end{array}\right)\right):
\; M\in{\cal D}_\delta,\; \vecxi\in\widetilde\Omega,\; \veca\in B_k(r) \right\} ,
$$
where ${\cal D}_\delta$ is an open subset of $\SL(2,\RR)$
specified below.

In the Iwasawa parametrization
3.8
$$
M=
\left( \begin{array}{cc}
uv^{-1/2} \sin\phi+v^{1/2}\cos\phi & uv^{-1/2} \cos\phi-v^{1/2}\sin\phi\\
v^{-1/2}\sin\phi& v^{-1/2}\cos\phi 
\end{array}\right) ,
$$
we have
\begin{eqnarray*}
 \Sigma_\delta(r)&=&
\Gamma 
\left\{\left(\tau,\phi;\vecxi +\left( \begin{array}{c} 
(u+v\cot\phi) \veca \\ 
 \veca \end{array}\right)\right):\right.  \\
&& \left.\phantom{\begin{array}{c} a\\b\end{array}}
\; (\tau,\phi)\in {\cal D}_\delta,\;  \vecxi\in\widetilde\Omega,\; 
\veca\in B_k(r v^{-1/2}\sin\phi ) \right\},
\end{eqnarray*}
where ${\cal D}_\delta$ is now chosen to be the open set of elements
$(\tau,\phi)\in\SL(2,\RR)$
subject to the restrictions
$$
0 < u + v\cot\phi < \delta, \quad
-1 < v^{-1/2} \cos\phi < 1 ,\quad
1 < v^{-1/2}\sin\phi < 2 .
$$

For the set
$$
\Pi_\delta(r)=
\Gamma 
\left\{\left(\tau,\phi;\vecxi +\left( \begin{array}{c} 
(u+v\cot\phi) \veca \\ 
\veca \end{array}\right)\right):  \\
\; (\tau,\phi)\in {\cal D}_\delta,\;
\vecxi\in\widetilde\Omega,\; \veca\in B_k(r) \right\},
$$
we find $\Pi_\delta(r)\subset \Sigma_\delta(r) \subset \Pi_\delta(2r)$.
Let us finally define
\begin{eqnarray*}
&&\hskip-22pt\widehat\Pi_{\varepsilon,\delta}(r)\\
&&\hskip-8pt=
\Gamma 
\left\{\left(\tau,\phi;\veczeta+\vecxi +\left( \begin{array}{c} 
\vecnull \\ 
\veca \end{array}\right)\right):  
  (\tau,\phi)\in {\cal D}_\delta,\;  
\veczeta\in B_{2k}(\varepsilon),\; \vecxi\in\widetilde\Omega,\; \veca\in B_k(r) \right\},
\end{eqnarray*}
where $B_{2k}(\varepsilon)\subset\RR^{2k}$ is the open ball of radius $\varepsilon$
about the origin. Thus,  $\widehat\Pi_{\varepsilon,\delta}(r)$ is   a full
dimensional (but thin) open set, 
which contains $\Pi_\delta(r)$ if $\delta>0$ is
chosen small enough. That is, for any $\varepsilon>0$ there is a $\delta>0$
such that
$$
\Sigma_\delta(r) \subset \widehat\Pi_{\varepsilon,\delta}(2r) .
$$

  5.11.9.
The characteristic function of $\widehat\Pi:=\widehat\Pi_{\varepsilon,\delta}(2r)$
therefore satisfies
$\chi_{\widehat\Pi}(\tau,\phi;\vecxi)=1$ for all
$(\tau,\phi;\vecxi)\in\Sigma_\delta(r)$. Hence, and because 
$\nu_H|_{\Gamma\backslash\Sigma(r)}\neq 0$,
there is a constant $c_\flat>0$ which is independent 
of $\delta$ and $\varepsilon$, such that, for all $\varepsilon>0$, $\delta>0$
sufficiently small,
$$
\nu_H(\Gamma\backslash\widehat\Pi)
=\int \chi_{\widehat\Pi}\;   d\nu_H \geq c_\flat 
\int_{{{0< u + v\cot\phi < \delta \atop
-1 < v^{-1/2} \cos\phi < 1} \atop
1 < v^{-1/2}\sin\phi < 2}}
\frac{du\,dv\,d\phi}{v^2} = 4c_\flat\, \delta
$$
and so  
$$
\nu(\Gamma\backslash\widehat\Pi)
\geq 
\nu_H(\Gamma\backslash\widehat\Pi)
\geq 
4c_\flat\, \delta.
$$
Since $\Gamma\backslash\widehat\Pi$ is open, we have 
along the subsequence $t_1,t_2,\ldots$ in 5.11.1
(Theorem~1, p.~311 in \cite{Shiryaev95})
$$
\liminf_{i\rightarrow\infty}\rho_{t_i}(\chi_{\widehat\Pi})
\geq 
\nu(\Gamma\backslash\widehat\Pi) 
\geq 
4c_\flat\, \delta.
$$

\vglue8pt 5.11.10.
On the other hand,  
$$
\chi_{\widehat\Pi}(\tau,\phi;\vecxi)
\leq F_{\varepsilon,\delta}(\tau,\phi;\vecxi)=
\sum_{\gamma\in\SL(2,\ZZ)} f_\delta(\gamma\tau,\phi_\gamma)
\eta_{\varepsilon}(\gamma\vecxi) ,
$$
where (as in~5.8)
$$
f_\delta(\tau,\phi)=\chi_1(u +v \cot\phi)\;
\chi_2(v^{-1/2}\cos\phi)\;
\chi_3(v^{-1/2}\sin\phi)
$$
and $\chi_1$, $\chi_2$, $\chi_3$ are the characteristic functions
of the intervals $[0,\delta]$, $[-1,1]$, $[1,2]$, respectively.
The function $\eta_{\varepsilon}$ is the characteristic
function of the set 
$$
\left\{\left(\veczeta+\vecxi +\left( \begin{array}{c} 
\vecnull \\ 
\veca \end{array}\right)\right): \; 
\veczeta\in B_{2k}(\varepsilon),\; 
\vecxi\in\widetilde\Omega,\; \veca\in B_k(r) \right\}
+\ZZ^{2k} .
$$

By Lemma~5.9, there is a constant $c_\sharp>0$
which is independent of $\delta$ and~$\varepsilon$, such that
$$
\limsup_{i\rightarrow \infty} \rho_{t_i}(\chi_{\widehat\Pi})
\leq c_\sharp \delta \int_{\TT^{2k}} \eta_{\varepsilon}(\vecxi)d\vecxi 
$$ 
for all sufficiently small $\varepsilon,\delta>0$.

We conclude that 
$$
4c_\flat \leq 
c_\sharp \int_{\TT^{2k}} \eta_{\varepsilon}(\vecxi)d\vecxi .
$$
This contradicts our assumption that $c_\flat>0$, if we can show that
the integral over $\eta_{\varepsilon}$ tends to zero,
as $\varepsilon\rightarrow 0$.
We will check this by a dimension consideration.

\vglue8pt 5.11.11.
To this end we need to show that, if $\Omega\neq\RR^{2k}$, we have
$$
\dim \left\{\left(\vecxi +\left( \begin{array}{c} 
\vecnull \\ 
\veca \end{array}\right)\right): \; 
\vecxi\in\widetilde\Omega,\; \veca\in B_k(r) \right\}<2k .
$$

In view of 5.11.7 this holds if and only if 
the dimension of the linear space
$$
V= \Omega + W , \qquad
W=\left\{\left( \begin{array}{c} 
\vecnull \\ 
\veca \end{array}\right): \; \veca\in\RR^k \right\},
$$
is strictly less than $2k$. Suppose $\dim V=\kappa$,  
and let $\vecb_1,\ldots,\vecb_\lambda$ form a basis of~$\Omega$.
Then there exist vectors $\vecb_{\lambda+1},\ldots,\vecb_\kappa\in W$
such that $\vecb_1,\ldots,\vecb_\kappa$ is a basis of $V$.
Hence 
$$
V=\Omega \oplus U, \quad U=\span\{\vecb_{\lambda+1},\ldots,\vecb_\kappa\}.
$$
The linear subspace 
$$
U^*=\left( \begin{array}{cc} 0 & -1 \\ 1 & 0 \end{array}\right) U \subset
\left( \begin{array}{c} \RR^k \\ \vecnull \end{array}\right) 
$$
clearly satisfies $U\cap U^*=\{\vecnull\}$, and also
$U^*\cap\Omega=\{\vecnull\}$ since $U\cap\Omega=\{\vecnull\}$
and $\Omega$ is $\SL(2,\RR)$-invariant.
Hence
$$
V \oplus U^* \subset \RR^{2k} ,
$$
and so $\dim V=2k$ implies $\dim U^*=\dim U=0$, which occurs only 
if $\Omega=V$.
Thus $\dim\Omega<2k$ implies $\dim V<2k$ and the claim is proved.

\vglue8pt 5.11.12.
Therefore $\nu_H\neq 0$ if and only if $H=G^k$, and hence the only
limit measure of converging subsequences is the normalized $\mu$.
The uniqueness of the limit measure implies finally that every subsequence
converges \cite{Shiryaev95}. 
\hfill\qed

\section{Diophantine conditions}\label{secdiophantine}

6.1. So far, all equidistribution results are valid only in the
case of bounded test functions $F$. We will now extend
these results to unbounded test functions~$F$, which grow moderately
in the cusps of $\Gamma\backslash G^k$. This will, however, only
be possible under certain diophantine assumptions on $\vecy$.

\vglue12pt 6.2. To this end let us discuss the following model situation.
Let $G=G^1$ and $\Gamma=\SL(2,\ZZ)\ltimes \ZZ^2$. Define
furthermore the subgroup
$$
\Gamma_\infty = \left\{ \left( \begin{array}{cc} 1 & m \\ 0 & 1 \end{array}\right):
m\in\ZZ\right\} \subset \SL(2,\ZZ),
$$
and put 
$$
v_\gamma := \Im (\gamma\tau) = \frac{v}{|c\tau+d|^2}, \quad
\hbox{ for } \gamma=\left( \begin{array}{cc} a & b \\ c & d \end{array}\right),
$$
and
$$
y_\gamma := \left( \begin{array}{c} 0 \\ 1 \end{array}\right)
\cdot(\gamma\vecxi) = c x + d y , \quad
\hbox{ with }\gamma\vecxi = \gamma\left( \begin{array}{c} x \\ y \end{array}\right)
= \left( \begin{array}{c} ax+by \\ cx+dy \end{array}\right).
$$
Let $\chi_R$ be the characteristic function of the interval $[R,\infty)$,
$$
\chi_R(t)=
\left\{ \begin{array}{ll} 
1 & (t\geq R) \\
0 & (t < R) .
\end{array}\right. 
$$
For any $f\in\C(\RR)$, which is rapidly decreasing at
$\pm\infty$, and $\beta\in\RR$, the function
$$
F_R(\tau;\vecxi)= 
\sum_{\gamma\in\Gamma_\infty\backslash\SL(2,\ZZ)}\, \sum_{m\in\ZZ} 
f\big( (y_\gamma+m) v_\gamma^{1/2} \big)\,
v_\gamma^\beta \, \chi_R(v_\gamma)
$$
is readily seen to be invariant under the action of $\Gamma$.
If $\tau$ lies in the fundamental domain of $\SL(2,\ZZ)$
given by 
${\cal F}_{\SL(2,\ZZ)}
=\{ \tau\in\H: u\in[-\tfrac12,\tfrac12), |\tau|> 1 \}$,
and if furthermore $R>1$, then $F_R(\tau;\vecxi)$ 
clearly has the representation
$$
F_R(\tau;\vecxi)= 
\sum_{m\in\ZZ} 
\big\{ f\big( (y+m) v^{1/2} \big) +f\big((-y+m) v^{1/2} \big) \big\}
v^\beta \chi_R(v) .
$$
The sum over $m$ is rapidly converging because $f$
is rapidly decreasing at $\pm\infty$.

We note that $F_R$ can alternatively be represented as
$$
F_R(\tau;\vecxi)=\sum_{(\gamma;\vecn)\in\widehat\Gamma_\infty\backslash\Gamma}
f\left( \left( \begin{array}{c} 0 \\ 1 \end{array}\right)
\cdot(\gamma\vecxi+\vecn) v_\gamma^{1/2}\right)\,
v_\gamma^\beta \, \chi_R(v_\gamma)
$$
with the abelian subgroup 
$$
\widehat\Gamma_\infty
=\left\{ \left(\left( \begin{array}{cc} 1 & m \\ 0 & 1 \end{array}\right);
 \left( \begin{array}{c} n \\ 0 \end{array}\right)\right) : \, m,n\in\ZZ\right\} \subset\Gamma.
$$

\vglue12pt 6.3.  We will assume from now on that $f\geq 0$.
The $\L^1$ norm of $F_R$ over $\Gamma\backslash G$ is then
$$
\mu(F_R)=\int_{\Gamma\backslash G} F_R(\tau;\vecxi)\; d\mu(\tau,\phi;\vecxi)
$$
with Haar measure 
$$
d\mu(\tau,\phi;\vecxi)=\frac{du\,dv\,d\phi\,dx\,dy}{v^2}.
$$
Then
$$
\mu(F_R)=\int_{\widehat\Gamma_\infty\backslash G} 
f\big( y  v^{1/2}\big)\,
v^\beta \, \chi_R(v) \; d\mu(\tau,\phi;\vecxi)
$$
and so
$$
\mu(F_R)=2\pi \int_\RR  f(w) dw  \int_R^\infty v^{\beta-5/2} dv
= 2\pi \frac{R^{-(3/2-\beta)}}{3/2-\beta}  \int_\RR  f(w) dw
$$
for $\beta<3/2$, and $\mu(F_R)=\infty$ otherwise.
Of special interest will be the case $\beta=1$, for which
$$
\mu(F_R)= 4\pi R^{-1/2} \int_\RR  f(w) dw .
$$

 6.4.  There is a well known one-to-one correspondence
between the coset $\Gamma_\infty\backslash\Gamma$ and the set
$$
\{ (0,1), (0,-1), (1,0), (-1,0) \} 
\cup \{ (c,d)\in\ZZ^2 :  c,d\neq 0, \; \gcd(c,d)=1 \} ,
$$
given by
$$
\left( \begin{array}{cc}
a & b \\
c & d 
\end{array}\right)
\mapsto
(c,d) .
$$
We may therefore write
\begin{eqnarray*}
&&\hskip-.65in
F_R(\tau;\vecxi) =
\sum_{m\in\ZZ} \big\{ f\big( (y+m) v^{1/2} \big) 
+ f\big( (-y+m) v^{1/2} \big) \big\}
v^\beta \chi_R(v) \\
&&\hskip-.3in+
\sum_{m\in\ZZ} 
\left\{ f\left((x + m) \frac{v^{1/2}}{|\tau|} \right) +
f\left((-x + m) \frac{v^{1/2}}{|\tau|} \right) \right\}
\frac{v^\beta}{|\tau|^{2\beta}} \chi_R\left(\frac{v}{|\tau|^2}\right) \\
&&\hskip-.3in+
\sum_{{{(c,d)\in\ZZ^2\atop \gcd(c,d)=1} \atop c,d\neq 0}} \,\sum_{m\in\ZZ} 
f\left((c x + d y + m) \frac{v^{1/2}}{|c\tau+d|}\right)
\frac{v^\beta}{|c\tau+d|^{2\beta}} \chi_R\left(\frac{v}{|c\tau+d|^2}\right).
\end{eqnarray*}

 From here on, we will only consider the case $\beta=1$,
and $\vecxi=\trans(0,y)$. 

\nonumproclaim{6.5. Proposition} Suppose $h\in\C(\RR)$ is positive and has compact support{\rm ,} and 
let $y$ be diophantine of type $\kappa$.
Then{\rm ,} for any $R>1$ and $\varepsilon,\varepsilon'$ with
$0<\varepsilon<1$ and $0<\varepsilon'<\frac{1}{\kappa-1}${\rm ,} 
$$
\limsup_{v\rightarrow 0}
\int_{|u|> v^{1-\varepsilon}} F_R\left(u+\i v; \left( \begin{array}{c} 0\\ y \end{array}\right)\right)
\; h(u)\, du
= O_{\varepsilon,\varepsilon'} (R^{-\varepsilon'/2}) ,
$$
where $y$, $f$ and $h$ are fixed.
\endproclaim

The proof of this proposition requires the following lemma.

\nonumproclaim{6.6. Lemma}  Let $\alpha$ be diophantine of type $\kappa$, and $f\in\C(\RR)$
be rapidly decreasing at $\pm\infty$ and positive{\rm ,} $f\geq 0$.
Then{\rm ,} for any fixed $A>1$ and $0<\varepsilon<\frac{1}{\kappa-1}${\rm ,}
$$
\sum_{d=1}^D \sum_{m\in\ZZ}  f\big( T (d\alpha+m) \big) 
\ll
\left\{ \begin{array}{ll}
T^{-A} & (D \leq T^\varepsilon) \\[10pt]
1 & (T^\varepsilon \leq D \leq T^{\frac{1}{\kappa-1}}) \\[10pt]
D\, T^{-\frac{1}{\kappa-1}}   & (D \geq T^{\frac{1}{\kappa-1}}) ,
\end{array}\right. 
$$
uniformly for all \pagebreak $D,T>1$.
\endproclaim

6.7. {\it Proof.}
\vglue8pt
6.7.1. 
Order $\alpha, 2\alpha, \ldots, D\alpha$ mod 1 in the unit interval $[0,1]$, 
and denote these
numbers by $0 < \varphi_1  < \ldots <\varphi_D < 1$.
Clearly $\varphi_{j+1}-\varphi_j = k_j\alpha\bmod 1$ for some 
integer $k_j\in[-D,D]$;
therefore, and because $\alpha$ is of type $\kappa$,
$$
\varphi_{j+1}-\varphi_j \geq \frac{C}{|k_j|^{\kappa-1}}\geq 
\frac{C}{D^{\kappa-1}} ,
$$
for some suitable constant $C>0$.
Hence in any interval of length $\ell$ there can be at most
$O(D^{\kappa-1}\ell+ 1)$ points.

\vglue8pt
6.7.2.  
As to the first bound,
take $\chi_{[-R,R]}$ to be the characteristic function of 
the interval $[-R,R]$ with $R>1$.
Then 
$$
\sum_{d=1}^D \sum_{m\in\ZZ} \chi_{[-R,R]}\big( T (d\alpha+m) \big) 
= 0 
$$
for $\frac{C T}{D^{\kappa-1}} > R$, since
$|d\alpha+m|\geq\frac{C}{d^{\kappa-1}}\geq\frac{C}{D^{\kappa-1}}$.
The argument in 6.7.1 shows that
$$
\sum_{d=1}^D \sum_{m\in\ZZ} \chi_{[-R,R]}\big( T (d\alpha+m) \big)
=O(R D^{\kappa-1}T^{-1} + 1 ) = O(R) 
$$
for $D^{\kappa-1}T^{-1}\leq 1$; hence
$$
\sum_{d=1}^D \sum_{m\in\ZZ} \chi_{[-R,R]}\big( T (d\alpha+m) \big)
=
\left\{ \begin{array}{ll}
O(R) & \hbox{ if $R \geq C \frac{T}{D^{\kappa-1}}$, } \\
0    & \hbox{ if $R < C \frac{T}{D^{\kappa-1}}$, } 
\end{array}\right. 
$$
in the range $D\leq T^{\frac{1}{\kappa-1}}$.
Since $f$ is rapidly decreasing, we have for any $B>3$
$$
f(t) \ll \sum_{R=1}^\infty R^{-B} \chi_{[-R,R]}(t),
$$ 
and hence, when $D\leq T^{\frac{1}{\kappa-1}}$,
$$
\sum_{d=1}^D \sum_{m\in\ZZ} f\big( T (d\alpha+m) \big)
\ll \sum_{R\geq C\frac{T}{D^{\kappa-1}}}^\infty R^{-(B-1)} \ll 
\left(\frac{D^{\kappa-1}}{T}\right)^{B-2}
$$ 
which proves the first bound in the range 
$D\leq T^\varepsilon\leq T^{\frac{1}{\kappa-1}}$,
for $\varepsilon<\frac{1}{\kappa-1}$.

\vglue8pt 
6.7.3.  
To prove the second and third relation, we follow \cite[pp.~13--14]{Davenport62}. Given any positive integer $q$
(to be fixed later) divide the sum over $d$ into blocks of the form
$$
\sum_{d=b}^{b+q-1} \sum_{m\in\ZZ} f\big( T (d\alpha+m) \big) 
=\sum_{d=0}^{q-1} \sum_{m\in\ZZ} f\big( T (b\alpha+d\alpha+m) \big) .
$$
(The last block might contain less than $q$ terms, but this is irrelevant
since we are seeking an upper bound.)
There are $O(\frac{D}{q} + 1)$ such blocks.
Take a rational approximation $\frac{p}{q}$ to $\alpha$ with 
$|\alpha-\frac{p}{q}|\leq q^{-2}$ and $p,q$ coprime, then the above sum is
$$
\sum_{d=0}^{q-1} \sum_{m\in\ZZ}  
f\left( T \left(b\alpha+\frac{dp+O(1)}{q}+m\right) \right) .
$$
Since $dp$ runs through a full set of residues mod $q$, the above
equals
$$
\sum_{r=0}^{q-1} \sum_{m\in\ZZ}  
f\left( T \left(b\alpha+\frac{r+O(1)}{q}+m\right)\right) 
=\sum_{r\in\ZZ} 
f\left( \frac{T}{q} (qb\alpha+r+O(1)) \right) .
$$
The term $qb\alpha$ may be replaced by the nearest integer
$+O(1)$, and so 
$$\sum_{r\in\ZZ} 
f\left( \frac{T}{q} (qb\alpha+r+O(1)) \right) 
=\sum_{r\in\ZZ} 
f\left( \frac{T}{q} (r+O(1)) \right)
$$
which in turn is clearly bounded by $O(\frac{q}{T} + 1)$
for $f$ is rapidly decreasing.
Therefore
$$
\sum_{d=1}^{D} \sum_{m\in\ZZ} f\big( T (d\alpha+m) \big) 
\ll \left(\frac{D}{q}+1\right)\left(\frac{q}{T}+1\right).
$$

\vglue8pt 
6.7.4.
By Dirichlet's theorem, we may take $\frac{p}{q}$ such that $q\leq T$ and 
$|\alpha-\frac{p}{q}|\leq q^{-1}T^{-1}$. Since $\alpha$ is of type $\kappa$,
we have $C q^{-\kappa} \leq |\alpha-\frac{p}{q}|$, so that
$$
T^{\frac{1}{\kappa-1}} \ll q \leq T ,
$$
and finally
\vglue12pt
\hfill ${\displaystyle
\sum_{d=1}^{D} \sum_{m\in\ZZ} f\big( T (d\alpha+m) \big) 
\ll \frac{D}{T^{\frac{1}{\kappa-1}}} +1 .
}$\hfill
\qed

\phantom{wanna}
 
\demo{{\rm 6.8.} Proof of Proposition~{\rm 6.5}}
\enddemo

6.8.1.
Because we are only concerned with upper bounds,
we may assume in the following without loss of generality
that $f$ is positive and even, i.e., $f\geq 0$, $f(-w)=f(w)$.

It follows from the expansion in 6.4 that, for $v<1$, 
the first term is absent, since $\chi_R(v)=0$ (recall: $R>1$); 
hence we are left with 
\begin{eqnarray*}
&&\hskip-.3in F_R\left(\tau;\left( \begin{array}{c} 0 \\ y\end{array}\right)\right)=
2 \sum_{m\in\ZZ} 
f\left(m \frac{v^{1/2}}{|\tau|} \right) 
\frac{v}{|\tau|^{2}} \chi_R\left(\frac{v}{|\tau|^2}\right) \\
&&\hskip.5in +2\sum_{{{(c,d)\in\ZZ^2\atop \gcd(c,d)=1}\atop c>0,d\neq 0}} \, \sum_{m\in\ZZ} 
f\left((d y + m) \frac{v^{1/2}}{|c\tau+d|} \right)
\frac{v}{|c\tau+d|^{2}} \chi_R\left(\frac{v}{|c\tau+d|^2}\right).
\end{eqnarray*} 

 6.8.2. 
As to the first term in the above expansion, a simple change of variable
$u=vt$ shows that
\begin{eqnarray*}
&&\int_{|u|>v^{1-\varepsilon}}   2 \sum_{m\in\ZZ} 
f\left(m \frac{v^{1/2}}{|\tau|} \right) 
\frac{v}{|\tau|^{2}} \chi_R\left(\frac{v}{|\tau|^2}\right) h(u)\, du \\
&&\qquad = 2 \int_{|t|>v^{-\varepsilon}} \sum_{m\in\ZZ} 
f\left(\frac{m}{v^{1/2} (t^2+1)^{1/2}} \right) 
\frac{1}{t^2+1} \chi_R\left(\frac{1}{v(t^2+1)}\right) h(vt)\, dt \\
&&\qquad \ll_{f,h} \int_{|t|>v^{-\varepsilon}} \frac{dt}{t^2+1} ,
\end{eqnarray*}
since the sum over $m$ is converging uniformly with respect to 
$t$ and $v$ due to the fact that $v(t^2+1)\leq R^{-1}< 1$.  

For $\varepsilon>0$
the value of the above integral converges to zero as $v\rightarrow 0$.

\vglue8pt 6.8.3.
We obtain an upper bound for the remaining terms, by dropping
the condition $|u|>v^{1-\varepsilon}$ in the integral. We are thus led
to estimate
$$
S(v)=
\sum_{{{(c,d)\in\ZZ^2\atop \gcd(c,d)=1}\atop c>0,d\neq 0}} \, \sum_{m\in\ZZ} 
J(v,c,d,m)
$$
with
$$
J(v,c,d,m)=
\int_\RR f\left((d y + m) \frac{v^{1/2}}{|c\tau+d|} \right)
\frac{v}{|c\tau+d|^{2}} \chi_R\left(\frac{v}{|c\tau+d|^2}\right) \; h(u)\, du.
$$
We substitute  $t=v^{-1}(u+\frac{d}{c})$ for $u$, yielding
$$
\frac{1}{c^2}
\int_\RR f\left((d y + m) \frac{1}{\sqrt{c^2 v(t^2+1)}} \right)
\frac{1}{t^2+1} \chi_R\left(\frac{1}{c^2v(t^2+1)}\right) \; h(vt-\frac{d}{c})\, dt.
$$
The range of integration is bounded by
$$
R < \frac{1}{c^2 v(t^2+1)}; \qquad\hbox{i.e.,}\quad
|t| \ll \frac{1}{c \sqrt{vR}}.
$$
This implies $|vt|\ll v^{1/2}c^{-1} R^{-1/2}$ is 
uniformly close to zero, and hence, because of the compact support of $h$,
we find $|d|\leq M c$, for some constant $M>0$ depending only
on the support of $h$. Therefore
$$
S(v) \ll
\sum_{c=1}^\infty\,
\sum_{0<|d|\leq M c}\;
\sum_{m\in\ZZ} 
K(v,c,d,m) ,
$$
with 
$$
K(v,c,d,m)=
\frac{1}{c^2}
\int_\RR f\left((d y + m) \frac{1}{\sqrt{c^2 v(t^2+1)}} \right)
\frac{1}{t^2+1} \chi_R\left(\frac{1}{c^2v(t^2+1)}\right) \; dt.
$$

 6.8.4.
In order to apply Lemma~6.6 with $D=M c$,
$T=(c^2 v(t^2+1))^{-1/2}>\sqrt{R}>1$, we split the $t$-range of integration
into the ranges
\begin{eqnarray*}
(1): &&   M c \leq (c^2 v(t^2+1))^{-\varepsilon'/2} \\
(2):&&   (c^2 v(t^2+1))^{-\varepsilon'/2} 
\leq M c \leq (c^2 v(t^2+1))^{-\frac{1}{2(\kappa-1)}} \\
(3): &&   M c  \geq (c^2 v(t^2+1))^{-\frac{1}{2(\kappa-1)}},
\end{eqnarray*}
which correspond to
\begin{eqnarray*}
(1): &&   D \leq T^{\varepsilon'} \\
(2): &&   T^{\varepsilon'} 
\leq D \leq T^{\frac{1}{\kappa-1}} \\
(3): &&   D  \geq T^{\frac{1}{\kappa-1}}.
\end{eqnarray*}
Here, $\varepsilon'<\frac{1}{\kappa-1}$.

We denote the corresponding integrals by
$K_1(v,c,d,m)$, $K_2(v,c,d,m)$ 
and $K_3(v,c,d,m)$, respectively.

\vglue8pt 6.8.5.
Because $R^{-1/2}\geq T^{-1}$, 
\begin{eqnarray*}
\sum_{c>0} \sum_{0<|d|\leq M c}\; \sum_{m\in\ZZ} 
K_1(v,c,d,m) 
& \dnhs\ll\dnhs&  R^{-A/2} \sum_{c>0}\frac{1}{c^2}
\int_{(1)}   \frac{1}{t^2+1} \chi_R\left(\frac{1}{c^2v(t^2+1)}\right) \;  dt \\
&\dnhs \ll\dnhs&  R^{-A/2} \sum_{c>0}\frac{1}{c^2}
\int   \frac{1}{t^2+1} \;  dt \\
&\dnhs \ll\dnhs& R^{-A/2} .
\end{eqnarray*}

\vglue8pt 6.8.6.
In order to obtain an upper bound,
we can relax the second range $T^{\varepsilon'}\leq D \leq T^\frac{1}{\kappa-1}$ 
to $R^{\varepsilon'/2}\leq D$, since $R^{1/2}\leq T$.
This yields
$$
\sum_{c>0}\, \sum_{0<|d|\leq M c}\; \sum_{m\in\ZZ} 
K_2(v,c,d,m) 
\ll  \sum_{M c\geq R^{\varepsilon'/2}} c^{-2}
\int   \frac{1}{t^2+1} \;  dt 
\ll R^{-\varepsilon'/2} .
$$

\vglue8pt 6.8.7.
In the third range, we find (putting $\delta=\frac{1}{\kappa-1}$),
\begin{eqnarray*}
&&\hskip-.75in \sum_{c>0} \, \sum_{0<|d|\leq M c}\; \sum_{m\in\ZZ} 
K_3(v,c,d,m)\\ & 
\ll  &\sum_{c>0} \frac{1}{c^2}
\int_{(3)} c^{1+\delta} v^{\delta/2} (t^2+1)^{\frac{\delta}{2}-1} 
\chi_R\left(\frac{1}{c^2v(t^2+1)}\right) \; dt \\ 
& \leq& \sum_{c>0} c^{-1+\delta} v^{\delta/2}
\int_\RR  
(t^2+1)^{\frac{\delta}{2}-1} 
\chi_R\left(\frac{1}{c^2v(t^2+1)}\right) \;  dt \\
& = &v^{\delta/2} \int_\RR \left\{ \sum_{c=1}^\infty c^{-1+\delta} 
\chi_R\left(\frac{1}{c^2v(t^2+1)}\right) \right\} (t^2+1)^{\frac{\delta}{2}-1} \, dt .
\end{eqnarray*}

For the inner sum there exist the upper bounds
\begin{eqnarray*}
&&\hskip-.25in \sum_{c=1}^\infty c^{-1+\delta} 
\chi_R\left(\frac{1}{c^2v(t^2+1)}\right)
\ll
\int_0^\infty x^{-1+\delta} 
\chi_R\left(\frac{1}{x^2v(t^2+1)}\right) dx \\
&&\quad =
[v(t^2+1)]^{-\delta/2} \int_0^\infty x^{-1+\delta} 
\chi_R\left(\frac{1}{x^2}\right) dx
=
[v(t^2+1)]^{-\delta/2} 
\left\{ \frac{x^{\delta}}{\delta} \right\}_0^{R^{-1/2}} ,
\end{eqnarray*}
and so
$$
\sum_{c>0} \, \sum_{d\ll c}\,  \sum_{m\in\ZZ} 
K_3(v,c,d,m)
\ll \frac{R^{-\delta/2}}{\delta}
\int_\RR (t^2+1)^{-1} \, dt  
= \frac{\pi R^{-\delta/2}}{\delta}   .
$$
The proof of Proposition~6.5 is complete.
\hfill\qed

\section{Equidistribution and unbounded test functions}\label{secequi}

7.1. 
Let us define the characteristic function on $\Gamma\backslash G^k$
(cf.~the proof of Proposition~5.4):
$$
X_R(\tau)= 
\sum_{\gamma\in\{\Gamma_\infty\cup(-1)\Gamma_\infty\}\backslash\SL(2,\ZZ)}
\chi_R(v_\gamma),
$$
where $\chi_R$ is the characteristic function of $[R,\infty)$.

 \vglue12pt 7.2.  
We shall consider functions on $\Gamma\backslash G^k$, which
grow moderately in the cusps. To be more precise,
we will require that, for some fixed constant $L>1$, the function {\it $F$
is dominated by $F_R$}; that is, for all sufficiently
large $R>1$, 
$$
|F(\tau,\phi;\vecxi)|X_R(\tau) \leq L+F_R(\tau;\vecxi)
$$
uniformly for all $(\tau,\phi;\vecxi)\in G^k$.
The function $F_R(\tau;\vecxi)$ is now viewed as a function
on $G^k$ (rather than $G^1$ as in Section \ref{secdiophantine}); 
that is, for 
$$\vecxi=\trans(x_1,\ldots,x_k,y_1,\ldots, y_k)$$ we put
$$
F_R(\tau;\vecxi)= 
\sum_{\gamma\in\Gamma_\infty\backslash\SL(2,\ZZ)} \sum_{m\in\ZZ} 
f\big( (y_{1,\gamma}+m) v_\gamma^{1/2} \big)\,
v_\gamma \, \chi_R(v_\gamma)
$$
which
is invariant under $\SL(2,\ZZ)\ltimes\ZZ^{2k}$ (Section \ref{secdiophantine})
and thus also under $\Gamma$.
Note that $F_R(\tau;\vecxi)$ is constant with respect to 
$x_2,\ldots,x_k$ and $y_2,\ldots,y_k$.
Again, $f\in\C(\RR)$ is rapidly decreasing at $\pm\infty$, positive and even.

\nonumproclaim{7.3. Theorem} Let $\Gamma$ be a subgroup of $\SL(2,\ZZ)\ltimes\ZZ^{2k}$ of finite
index. Let $h$ be a continuous probability density $\RR\rightarrow\RR_+$
with compact support. 
Suppose the continuous function $F\geq 0$ is dominated by $F_R$.
Fix some $\vecy\in\TT^{k}$
such that the components of the vector 
$(\trans\vecy,1)\in\RR^{k+1}$ are linearly
independent over $\QQ$.
Then{\rm ,} for any $\varepsilon$ with $0<\varepsilon<1${\rm ,}
$$
\liminf_{v\rightarrow 0}
\int_{|u|> v^{1-\varepsilon}} F\left(u+\i v,0;
\left( \begin{array}{c}\vecnull \\ \vecy \end{array}\right)\right) \; h(u)\, du
\geq \frac{1}{\mu(\Gamma\backslash G^k)}\int_{\Gamma\backslash G^k} F \, d\mu .
$$
Assume furthermore
that $y_1$ is diophantine.
Then{\rm ,} for any $\varepsilon$ with $0<\varepsilon<1${\rm ,}
$$
\limsup_{v\rightarrow 0}
\int_{|u|> v^{1-\varepsilon}} F\left(u+\i v,0;
\left( \begin{array}{c}\vecnull \\ \vecy \end{array}\right)\right) \; h(u)\, du
\leq \frac{1}{\mu(\Gamma\backslash G^k)}\int_{\Gamma\backslash G^k} F \, d\mu .
$$
\endproclaim

\demo{Proof}
We obtain the {\it lower bound} from the function
$$
G_R(\tau,\phi; \vecxi):=F(\tau,\phi; \vecxi)\, (1-X_R(\tau))
\leq F(\tau,\phi; \vecxi).
$$
Clearly, $G_R$ is bounded. Therefore
$$
\int_{|u|> v^{1-\varepsilon}} G_R(u+\i v,0; \vecxi) \; h(u)\, du
=\int_\RR G_R(u+\i v,0; \vecxi) \; h(u)\, du +O_R(v^{1-\varepsilon}),
$$
and, by Theorem~5.7,\footnote{\label{foot} The fact that 
$G_R$ is only piecewise continuous should not worry us:
Theorem~5.7 can easily be extended to such functions
by approximating these
from above and from below by continuous functions. 
In any case, the argument presented
here works as well if $\chi_R$ is smoothed slightly, which makes
$G_R$ continuous.}
$$
\lim_{v\rightarrow0}
\int_\RR G_R(u+\i v,0; \vecxi) \; h(u)\, du
=\frac{1}{\mu(\Gamma\backslash G^k)}
\int_{\Gamma\backslash G^k} G_R \, d\mu .
$$
Now since $0\leq F X_R \leq LX_R + F_R$ for $R$ large enough
we have
$$
\int_{\Gamma\backslash G^k} F X_R\, d\mu 
\leq \int_{\Gamma\backslash G^k} (LX_R+F_R)\,  d\mu \ll L R^{-1}+R^{-1/2}
$$
from~5.4 and 6.3, and hence
$$
\int_{\Gamma\backslash G^k} G_R \, d\mu
=\int_{\Gamma\backslash G^k} F \, d\mu + O(LR^{-1}+R^{-1/2}).
$$
In summary
$$
\liminf_{v\rightarrow 0}
\int_{|u|> v^{1-\varepsilon}} F(u+\i v,0; \vecxi) \; h(u)\, du
\geq \frac{1}{\mu(\Gamma\backslash G^k)}\int_{\Gamma\backslash G^k} F \, d\mu 
+O(R^{-1/2}),
$$
for all $R$ large enough. The assertion on the lower bound
follows now from the fact that $R$ can be chosen arbitrarily large.

For the {\it upper bound}, notice that for $R$ large enough,
$$
F(\tau,\phi; \vecxi) \leq F(\tau,\phi; \vecxi) 
(1-X_R(\tau)) + LX_R(\tau)+ F_R(\tau; \vecxi) .
$$
By virtue of the bound obtained in the previous paragraph,
and by Proposition~6.5, we find
that
\begin{eqnarray*}
&&\hskip-.5in
\limsup_{v\rightarrow 0}
\int_{|u|> v^{1-\varepsilon}} F(u+\i v,0; \vecxi) \; h(u)\, du 
\\
&&\leq \frac{1}{\mu(\Gamma\backslash G^k)}\int_{\Gamma\backslash G^k} F \, d\mu 
+ O(R^{-1/2})+ O(R^{-\eta})
\end{eqnarray*}
for some small constant $\eta>0$. 
This holds again for arbitrarily large $R$,
and the statement is proved.
\enddemo

\nonumproclaim{7.4. {C}orollary} 
 Let $\Gamma$, $h$, $\vecy$ be as in Theorem~{\rm 7.3,}
and $F: \Gamma\backslash G^k\rightarrow\CC$ be
a continuous function which is dominated by $F_R$.
If $y_1$ is diophantine{\rm ,} then{\rm ,} for any $\varepsilon$ with $0<\varepsilon<1${\rm ,}
$$
\lim_{v\rightarrow 0}
\int_{|u|> v^{1-\varepsilon}} F\left(u+\i v,0;
\left( \begin{array}{c}\vecnull \\ \vecy \end{array}\right)\right) \; h(u)\, du
= \frac{1}{\mu(\Gamma\backslash G^k)}\int_{\Gamma\backslash G^k} F \, d\mu .
$$
\endproclaim

\demo{Proof}
Define
$$
\Re_+ F(\tau,\phi;\vecxi)
=
\left\{ \begin{array}{ll}
\Re F(\tau,\phi;\vecxi) & \hbox{if $\Re F(\tau,\phi;\vecxi)>0$,}\\
0 & \hbox{if $\Re F(\tau,\phi;\vecxi)\leq 0$,}
\end{array}\right. 
$$
and $\Re_- F=\Re_+ F-\Re F$. We similarly define $\Im_\pm F$
as the positive/negative part of $\Im F$.
Then 
$$
F = \Re_+ F -\Re_- F + \i \Im_+ F -\i \Im_- F
$$
with 
$$
0\leq \Re_+ F X_R\leq L+F_R,\quad
0\leq \Re_- F X_R\leq L+F_R,
$$
$$
0\leq \Im_+ F X_R\leq L+F_R,\quad
0\leq \Im_- F X_R\leq L+F_R .
$$
We can thus apply Theorem~7.3 to each term separately,
$$
\lim_{v\rightarrow 0}
\int_{|u|> v^{1-\varepsilon}} \Re_\pm F\left(u+\i v,0;
\left( \begin{array}{c}\vecnull \\ \vecy \end{array}\right)\right) \; h(u)\, du
= \frac{1}{\mu(\Gamma\backslash G^k)}\int_{\Gamma\backslash G^k} 
\Re_\pm F \, d\mu 
$$
and likewise for $\Im_\pm F$.
\enddemo

7.5.
Since in our main application $\Gamma=\Gamma^k$, which is a subgroup
of finite index in $\SL(2,\ZZ)\ltimes (\tfrac12\ZZ)^{2k}$
rather than in $\SL(2,\ZZ)\ltimes\ZZ^{2k}$ (Lemma~4.12), we restate Corollary
7.4 in the following equivalent way.
Define the dominating function $\hat F_R$ on 
$\Gamma\backslash G^k$ by 
$\hat F_R(\tau;\vecxi)=F_R(\tau; 2 \vecxi)$,
with $F_R$ as in 7.2.

\nonumproclaim{7.6. {C}orollary} 
 Let $\Gamma$ be a subgroup of $\SL(2,\ZZ)\ltimes (\tfrac12\ZZ)^{2k}$
of finite index{\rm ,} 
$h${\rm ,} $\vecy$ be as in Theorem~{\rm 7.3,}
and $F: \Gamma\backslash G^k\rightarrow\CC$
a continuous function which is dominated by $\hat F_R$.
If $y_1$ is diophantine{\rm ,} then{\rm ,} for any $\varepsilon$ with $0<\varepsilon<1${\rm ,}
$$
\lim_{v\rightarrow 0}
\int_{|u|> v^{1-\varepsilon}} F\left(u+\i v,0;
\left( \begin{array}{c}\vecnull \\ \vecy \end{array}\right)\right) \; h(u)\, du
= \frac{1}{\mu(\Gamma\backslash G^k)}\int_{\Gamma\backslash G^k} F \, d\mu .
$$
\endproclaim

\demo{Proof}
Apply Corollary 7.4 with the test function
$\tilde F:\tilde\Gamma\backslash G^k\rightarrow \CC$
defined by
$$
\tilde F(\tau,\phi;\vecxi)=F(\tau,\phi;\tfrac12 \vecxi)
$$
where 
$$
\tilde\Gamma=
\left(\left( \begin{array}{cc} 
2 & 0 \\ 0 & 2
\end{array}\right); \vecnull\right)\;
\Gamma\;
\left(\left( \begin{array}{cc} 
\tfrac12 & 0 \\ 0 & \tfrac12
\end{array}\right); \vecnull\right)
$$
is a subgroup of finite index in $\SL(2,\ZZ)\ltimes\ZZ^{2k}$
(compare Remark~4.13).
\enddemo

\section{The main theorem} \label{secmain}

\nonumproclaim{8.1. Main Theorem} Suppose $f(w_1,w_2)=\psi_1(w_1^2+w_2^2)$ and $g(w_1,w_2)
=\psi_2(w_1^2+w_2^2)$ with $\psi_1,\psi_2\in\Sw(\RR_+)$.
Let $h$ be a continuous function $\RR\rightarrow\CC$
with compact support. 
Assume that $y_1,y_2,1$ are linearly independent over $\QQ$ and that
$y_1$ is diophantine.
Then{\rm ,} with $\vecxi=\trans(0,0,y_1,y_2)${\rm ,}  
\begin{eqnarray*}
&&\hskip-.5in \lim_{v\rightarrow 0}
\int_\RR \Theta_f(u+\i v,0;\vecxi)
\overline{\Theta_g(u+\i v,0;\vecxi)} \; h(u)\, du \\
&&\qquad = \pi \big\{ 2\pi h(0) +\int_\RR h(u)\, du\big\} 
\int_0^\infty \psi_1(r)\overline{\psi_2(r)}\, dr . 
\end{eqnarray*}
\endproclaim
 
The proof of the main theorem requires the following two lemmas.

\nonumproclaim{8.2. Lemma}  If $f,g\in\Sw(\RR^2)${\rm ,}
$$
\frac{1}{\mu(\Gamma^2\backslash G^2)}
\int_{\Gamma^2\backslash G^2} 
\Theta_f(\tau,\phi;\vecxi)\overline{\Theta_g(\tau,\phi;\vecxi)}\, d\mu 
= \iint f(w_1,w_2)\overline{g(w_1,w_2)} \,dw_1\, dw_2 .
$$
\endproclaim

Note that if $f(w_1,w_2)=\psi_1(w_1^2+w_2^2)$ and 
$g(w_1,w_2)=\psi_2(w_1^2+w_2^2)$, then
$$
\iint f(w_1,w_2)\overline{g(w_1,w_2)} \,dw_1\, dw_2
= \pi \int_0^\infty \psi_1(r)\overline{\psi_2(r)}\, dr .
$$ 

\demo{Proof}
A short calculation shows that
$$
\int_{\TT^4}
\Theta_f(\tau,\phi;\vecxi)\overline{\Theta_g(\tau,\phi;\vecxi)} \, d\vecxi  
= \iint f_\phi(w_1,w_2)\overline{g_\phi(w_1,w_2)} \,dw_1\, dw_2 .
$$
Since $f_\phi=\tilde R(\i,\phi) f$ with $\tilde R(\i,\phi)$ unitary, we have
\vglue12pt
\hfill ${\displaystyle
\iint f_\phi(w_1,w_2)\overline{g_\phi(w_1,w_2)}\,dw_1\, dw_2
=\iint f(w_1,w_2)\overline{g(w_1,w_2)}\,dw_1\, dw_2 . 
}$
\enddemo
 \pagebreak

\nonumproclaim{8.3. Lemma}  Suppose $f(w_1,w_2)=\psi_1(w_1^2+w_2^2)$ and 
$g(w_1,w_2)=\break \psi_2(w_1^2+w_2^2)${\rm ,} with $\psi_1,\psi_2\in\Sw(\RR_+)$.
For any $\tfrac12<\gamma<1${\rm ,} 
\begin{eqnarray*}
&&\hskip-.51in
\lim_{v\rightarrow 0} \int_{|u|< v^{\gamma}}
\Theta_f\left(u+\i v, 0 ;\left( \begin{array}{c}\vecnull \\ \vecy \end{array}\right)\right)
\overline{\Theta_g\left(u+\i v, 0 ;
\left( \begin{array}{c}\vecnull \\ \vecy \end{array}\right)\right)} h(u)\,du  \\[4pt]
&&\qquad\qquad 
= 2 \pi^2 h(0) \int_0^\infty \psi_1(r)\overline{\psi_2(r)}\, dr .
\end{eqnarray*}
\endproclaim

\demo{Proof}
Proposition~4.11 tells us that
\begin{eqnarray*}
&&\hskip-.5in
\Theta_f\left(-\frac{1}{\tau},\arg\tau;
\left( \begin{array}{c} -\vecy \\ \vecnull \end{array}\right)\right)
\overline{\Theta_g\left(-\frac{1}{\tau},\arg\tau;
\left( \begin{array}{c} -\vecy \\ \vecnull \end{array}\right)\right)} 
\\[4pt]
&&\qquad = \frac{v}{|\tau|^2} f_{\arg\tau}(0,0) \overline{g_{\arg\tau}(0,0)}
+O_R\left(\left(\frac{v}{|\tau|^2}\right)^{-R}\right)
\end{eqnarray*}
holds uniformly for $\Im(-\tau^{-1})=v|\tau|^{-2}>\frac12$.
This condition is met, e.g., when $|u|<v^{1/2}<1$. 
For $|u|<v^{\gamma}<1$, with $\frac12<\gamma<1$, the error term is bounded by
$$
O_R\left(\left(\frac{v}{|\tau|^2}\right)^{-R}\right)=
O_R(v^{R(2\gamma-1)}).
$$
Now replacing $(w_1,w_2)$ by polar coordinates $(r\cos\zeta,r\sin\zeta)$
yields
\begin{eqnarray*}
f_{\arg\tau}(0,0) \overline{g_{\arg\tau}(0,0)} 
& =& \frac{|\tau|^2}{v^2}
\left\{\iint e\left(\tfrac12 (w_1^2+w_2^2)\frac{u}{v}\right)\, f(w_1,w_2)\, dw_1\, dw_2\right\}
  \\
&&  \times \overline{
\left\{\iint e\left(\tfrac12 (w_1^2+w_2^2)\frac{u}{v}\right)\, g(w_1,w_2)\, dw_1\, dw_2\right\}}\\
& =& \frac{|\tau|^2}{v^2} \pi^2 
\iint_0^\infty e\left(\frac{(r_1-r_2)u}{2v}\right) \psi_1(r_1)
\overline{\psi_2(r_2)}\, dr_1 dr_2\\
&=&\frac{|\tau|^2}{v^2} \pi^2 \hat\psi_1\left(\frac{u}{2v}\right)
\overline{\hat\psi_2\left(\frac{u}{2v}\right)},
\end{eqnarray*}
where $\hat\psi$ denotes the Fourier transform 
$$
\hat\psi(u)=\int_0^\infty e(ur) \psi(r)\, dr .
$$
Clearly $\hat\psi\in\L^2(\RR)$ for $\psi\in\Sw(\RR_+)\subset\L^2(\RR_+)$. 
Thus,
\begin{eqnarray*}
&&\hskip-.75in\int_{|u|< v^{\gamma}}  
\Theta_f\left(u+\i v,0;\left( \begin{array}{c}\vecnull \\ \vecy \end{array}\right)\right) 
\overline{\Theta_g\left(u+\i v,0;
\left( \begin{array}{c}\vecnull \\ \vecy \end{array}\right)\right)} h(u)\,du \\
& = &\frac{\pi^2}{v} \int_{|u|< v^{\gamma}}
\hat\psi_1\left(\frac{u}{2v}\right)\overline{\hat\psi_2\left(\frac{u}{2v}\right)}
 h(u)\,du + O_R(v^{\gamma+R(2\gamma-1)}) \\
& =& 2\pi^2 \int_{2|u|< v^{\gamma-1}}
\hat\psi_1(u)\overline{\hat\psi_2(u)}
 h(2vu)\,du + O_R(v^{\gamma+R(2\gamma-1)}).
\end{eqnarray*}
Since $h$ is continuous, for any given $\varepsilon>0$ we find
a $v_0>0$, such that
$$
|h(2vu)-h(0)|<\varepsilon, \hbox{ uniformly for all $2|u|< v^{\gamma-1}$,
$0<v<v_0$.}
$$
Thus for any $\varepsilon>0$
\begin{eqnarray*}
&&\hskip-.75in \lim_{v\rightarrow 0}\int_{|u|< v^{\gamma}}  
\Theta_f\left(u+\i v,0;\left( \begin{array}{c}\vecnull \\ \vecy \end{array}\right)\right) 
\overline{\Theta_g\left(u+\i v,0;
\left( \begin{array}{c}\vecnull \\ \vecy \end{array}\right)\right)} h(u)\,du \\
& = &\lim_{v\rightarrow 0}
2\pi^2 \{h(0)+O(\varepsilon)\}\int_{2|u|< v^{\gamma-1}}
\hat\psi_1(u)\overline{\hat\psi_2(u)} \,du \\
& = &
2\pi^2 \{h(0)+O(\varepsilon)\}\int_{\RR}
\hat\psi_1(u)\overline{\hat{\psi}_2(u)} \,du \\
& = &
2\pi^2 \{h(0)+O(\varepsilon)\}\int_0^\infty
\psi_1(r)\overline{\psi_2(r)} \,dr 
\end{eqnarray*}
by Parseval's equality. Because $\varepsilon>0$ can be arbitrarily small, 
the claim is proved.
\enddemo

 8.4. {\it Proof of the main theorem}.

\vglue8pt 8.4.1. 
Due to the linearity in $h$ of the integrals in 8.1, 
we may assume without 
loss of generality that (i) $h$ is positive 
(compare the argument used in the proof of Corollary 7.4)
and (ii) that $h$ is normalized as a probability density.

\vglue8pt 8.4.2.
Let us split the integration on the left-hand side of
8.1 into 
$$
\int_\RR = \int_{|u|< v^{1-\varepsilon}} + \int_{|u|> v^{1-\varepsilon}} ,
$$
for some small $\varepsilon>0$.
The first integral gives, by virtue of Lemma~8.3,
the contribution
$$
2\pi^2 h(0) \int_0^\infty \psi_1(r)\overline{\psi_2(r)}\, dr.
$$

\vglue8pt 8.4.3.
In order to apply Corollary~7.6, we need to construct a
function $F_R$ of the form studied in 7.2,
which dominates $|F|$.
Let us define
$$
f^*(w_1)=\sup_{w_2\in\RR}\sup_{\phi\in\RR} |f_\phi(w_1,w_2)\, g_\phi(w_1,w_2)|,
$$
which is clearly rapidly decreasing at $\pm\infty$ since, for every $T>1$,
there is a constant $c_T>0$ such that
$$
f^*(w_1)\leq
\sup_{w_2\in\RR} c_T \left(1+\sqrt{w_1^2+w_2^2}\right)^{-2T}
\leq c_T (1+|w_1|)^{-2T},
$$
holds (cf.~Lemma 4.3).

Choosing (compare 7.2)
$$
F_R(\tau;\vecxi)
=\sum_{\gamma\in\Gamma_\infty\backslash\SL(2,\ZZ)}\, \sum_{m\in\ZZ} 
f^* \big(-\tfrac12(y_{1,\gamma}+m) v_\gamma^{1/2} \big)\,
v_\gamma \, \chi_R(v_\gamma),
$$
we have for all $v>R$ 
$$
\hat F_R(\tau;\vecxi)=F_R(\tau;2\vecxi)=
v \sum_{m\in\ZZ} \left\{ f^* \left( \left(\tfrac{m}{2}-y_{1}\right) v^{1/2} \right)
+f^* \left( \left(\tfrac{m}{2}+y_{1}\right) v^{1/2}\right) \right\};
$$
that is,
$$
\hat F_R(\tau;\vecxi)=
v \left\{ f^* \left( \left(\tfrac{n}{2}-y_{1}\right) v^{1/2} \right)
+f^* \left( \left(-\tfrac{n}{2}+y_{1}\right) v^{1/2} \right)\right\} + O(v^{-T}),
$$
for all $y_1\in \tfrac{n}{2}+[-\tfrac14,\tfrac14]$, $n\in\ZZ$.
By construction, for $\vecn=\trans(n_1,n_2)$,
$$
\big|f_\phi\left( (\vecn-\vecy) v^{1/2} \right)
g_\phi\left( (\vecn-\vecy) v^{1/2} \right)\big|
\leq f^* \left( (n_1-y_{1}) v^{1/2}\right)
$$
which implies that, for all $v>R$, $R$ large enough,
\begin{eqnarray*}
&&\hskip-.25in v\Bigg|\sum_{\vecm\in\ZZ^2}
f_\phi\left((\vecm-\vecy) v^{1/2}\right)
g_\phi\left((\vecm-\vecy) v^{1/2} \right)\Bigg|\\
&&\quad =v\big|f_\phi\left( (\vecn-\vecy) v^{1/2} \right)
g_\phi\left( (\vecn-\vecy) v^{1/2} \right)\big|
+O(v^{-T})
\\
&&\quad \leq  vf^* \big( (n_1-y_{1}) v^{1/2} \big) +O(v^{-T})
=v\sum_{m_1\in\ZZ} f^* \big( (m_1-y_{1}) v^{1/2} \big)+O(v^{-T}) ,
\end{eqnarray*}
uniformly for $\vecy=\trans(y_1,y_2)\in\vecn+[\frac12,\frac12]^2$, 
$\vecn\in\ZZ^2$.
 
Therefore, by virtue of Proposition~4.10, we have,
for all sufficiently large~$R$,
$$
|\Theta_f(\tau,\phi;\vecxi)\Theta_g(\tau,\phi;\vecxi)| \leq 
1+ v \sum_{{m\in\ZZ\atop  m \ {\rm even}}} 
f^* \left( \left(\tfrac{m}{2}-y_{1}\right) v^{1/2} \right)
\leq 1+ \hat F_R(\tau;\vecxi) 
$$
for $v\geq R$, and so
$|\Theta_f \Theta_g|X_R \leq 1+\hat F_R$.
We can now apply Corollary~7.6, and thus obtain 
the second term on the right-hand side of 8.1 
(recall Lemma~8.2). \phantom{mountains}
\hfill\qed

\vglue12pt 8.5. {\it Proof of Theorem}~2.2.
Recall that
\begin{eqnarray*}
&&\hskip-.5pt
\int_\RR 
\Theta_f\left(u+\i \frac{1}{\lambda},0;\trans(0,0,\alpha,\beta) \right)
\overline{\Theta_g\left(u+\i \frac{1}{\lambda},0;\trans(0,0,\alpha,\beta) \right)} 
\; h(u)\, du\\[4pt]
&&\qquad\quad
= \pi R_2(\psi_1,\psi_2,h,\lambda) .
\end{eqnarray*}
We have furthermore 
$$
\hat h(s)= \int_\RR h(u) e\left(\tfrac12 us\right)\, du,
\quad 
h(u)= \tfrac12 \int_\RR \hat h(s) e\left(-\tfrac12 us \right)\, ds;
$$
hence $2h(0)=\int \hat h(s)ds$ and $\int h(u) du = \hat h(0)$.
\hfill\qed
\vglue12pt
8.6. {\it Proof of Theorem} 1.8.

\vglue8pt 8.6.1.
Let $\chi[a,b]$ be the characteristic function of the interval
$[a,b]$. Given any $\varepsilon>0$, we approximate $\chi[a,b]$
from above and below by 
functions $\chi_\pm\in\C^\infty(\RR)$ with compact support so that
$$
\chi_-(s) \leq \chi[a,b](s) \leq \chi_+(s), \qquad
\int_\RR (\chi_+(s) -\chi_-(s)) \, ds <\varepsilon.
$$
Put 
$$
\hat h_\pm(s) = \chi_\pm(s)  \pm \frac{\delta}{1+s^2} ,
$$
where $\delta>0$ is chosen such that
$$
\int_\RR \frac{4\delta}{1+s^2}\, ds < \varepsilon.
$$
Then 
$$
\hat h_-(s) +\frac{\delta}{1+s^2} \leq \chi[a,b](s) \leq 
\hat h_+(s) -\frac{\delta}{1+s^2},$$
$$
\int_\RR \left(\hat h_+(s) -\hat h_-(s)+\frac{2\delta}{1+s^2} \right)ds <2\varepsilon.
$$
The inverse Fourier transform
$$
h_\pm(u)= \tfrac12\int_\RR \hat h_\pm(s) e\left(-\tfrac12 us\right)\, ds
$$
is continuous on $\RR$, infinitely differentiable
on $\RR-\{0\}$ and decreases, together with its derivatives,
rapidly at $\pm\infty$.

\vglue8pt 8.6.2. 
We fix a smoothed characteristic function
$\chi\in\C^\infty(\RR)$ of compact support in $[-2,2]$, 
with $0\leq \chi\leq 1$ and $\chi(u)=1$ if $u\in[-1,1]$.
Define
$$
h_{T,\pm}(u)=h_\pm(u)\, \chi\left(\frac{u}{T} \right),
$$
which is continuous and has compact support in $[-2T,2T]$.
For the Fourier transform
$$
\hat h_{T,\pm}(s)= \int_\RR h_{T,\pm}(u) e\left(\tfrac12 us \right)\, du
$$
we have, for some constant $C$,
$$
|\hat h_\pm(s)-\hat h_{T,\pm}(s)| 
\leq 
\int_\RR |h_\pm(u)|\, \left|1-\chi\left(\frac{u}{T} \right)\right|\, du 
\leq
\int_{|u|>T} |h_\pm(u)|\, du
\leq 
\frac{C}{T},
$$
and (integrate by parts twice)
\begin{eqnarray*}
  |\hat h_\pm(s)-\hat h_{T,\pm}(s)|  
&\leq &
\frac{1}{(\pi s)^2}
\left\{\int_{|u|>T} |h_\pm''(u)|\, du 
+\frac2T \int_\RR \left|h_\pm'(u) \chi'\left(\frac{u}{T}\right)\right|\, du\right. \\[4pt]
&&\left.\phantom{\frac{1}{(\pi s)^2}
\{}
+\frac{1}{T^2}\int_\RR \left|h_\pm(u) \chi''\left(\frac{u}{T}\right)\right|\, du  \right\} \\
&
\leq&
\frac{C}{T s^2}.
\end{eqnarray*}
Therefore we find some $T>1$ such that
$$
|\hat h_\pm(s)-\hat h_{T,\pm}(s)| < \frac{\delta}{1+s^2} .
$$
Hence
$$
\hat h_{T,-}(s) \leq \chi[a,b](s) \leq \hat h_{T,+}(s),
\qquad
\int_\RR (\hat h_{T,+}(s) -\hat h_{T,-}(s)) ds <2\varepsilon.
$$

 8.6.3.
We will assume in the following that $\psi_1,\psi_2\geq 0$.
Then
\begin{eqnarray*}
&&\hskip-.75in\frac{1}{\pi\lambda} \sum_{j\neq k} 
\psi_1\left(\frac{\lambda_j}{\lambda} \right)\psi_2\left(\frac{\lambda_k}{\lambda} \right)
\hat h_{T,-}(\lambda_j-\lambda_k) 
\\
&&\leq
\frac{1}{\pi\lambda} \sum_{j\neq k} 
\psi_1\left(\frac{\lambda_j}{\lambda} \right)\psi_2\left(\frac{\lambda_k}{\lambda} \right)
\chi[a,b](\lambda_j-\lambda_k) \\ 
&&\leq
\frac{1}{\pi\lambda} \sum_{j\neq k} 
\psi_1\left(\frac{\lambda_j}{\lambda} \right)\psi_2\left(\frac{\lambda_k}{\lambda} \right)
\hat h_{T,+}(\lambda_j-\lambda_k) .
\end{eqnarray*}
The functions $h_{T,\pm}$ satisfy the assumptions in 
Theorem~2.2, so
the limits of left- and right-hand sides exist,
and differ by less than
$$
2\pi\varepsilon \; \big|\int_0^\infty \psi_1(r)\psi_2(r)dr\big|
$$
for arbitrarily small $\varepsilon>0$.
Hence 
$$
\lim_{\lambda\rightarrow\infty}
\frac{1}{\pi\lambda} \sum_{j\neq k} 
\psi_1\left(\frac{\lambda_j}{\lambda} \right)\psi_2\left(\frac{\lambda_k}{\lambda} \right)
\chi[a,b](\lambda_j-\lambda_k) =
\pi (b-a) \int \psi_1(r)\psi_2(r)dr.
$$

\vglue8pt 8.6.4.
Analogous arguments allow us to replace first $\psi_1$ 
and then $\psi_2$ by characteristic functions.
\hfill\qed
\vglue12pt
For detailed discussions of approximation functions 
of the type used above, see \cite{Vaaler85} and references therein.

\section{Counterexamples} \label{counter}

9.1. 
Put 
$$
Q_{\alpha,\beta}(m,n)= (m-\alpha)^2+(n-\beta)^2. 
$$
For $(\alpha,\beta)\in\QQ^2$ we find (see Appendix~A.10 for details)
for $\lambda\rightarrow\infty$,
$$
R^{(\alpha,\beta)}[0,0]=
\frac{1}{\pi \lambda} 
\#\{ (m_1,m_2,n_1,n_2)\in\ZZ^4 : (m_1,n_1)\neq (m_2,n_2), $$
$$Q_{\alpha,\beta}(m_1,n_1)  \leq  \lambda,\;
Q_{\alpha,\beta}(m_1,n_1)=Q_{\alpha,\beta}(m_2,n_2) \} 
\sim c_{\alpha,\beta} \log\lambda ,
$$
for some constant $c_{\alpha,\beta}>0$. This fact will be the key in
 proving the first half of Theorem~1.13.

\vglue12pt 9.2. {\it  Proof of Theorem}~1.13 (i).
Enumerate the rational forms $Q_{\alpha_j,\beta_j}$ with
$(\alpha_j,\beta_j)\in\QQ^2$ as $P_1,P_2,P_3,\ldots\, $. Because of the
asymptotics 9.1, given any $\lambda>1$, there exists
an $M_j>\lambda$ such that 
\begin{eqnarray*}
&&\frac{1}{\pi M_j} 
\#\{ (m_1,m_2,n_1,n_2)\in\ZZ^4 : (m_1,n_1)\neq (m_2,n_2), \\
&&\phantom{\frac{1}{\pi M_j} 
\#\{}P_j(m_1,n_1) \leq M_j ,\;
P_j(m_1,n_1)=P_j(m_2,n_2) \} 
\geq \frac{\log M_j}{\log\log\log M_j} .
\end{eqnarray*}
Now since
\begin{eqnarray*}
Q_{\alpha,\beta}(m_1,n_1)-Q_{\alpha,\beta}(m_2,n_2)
&\dnhs=\dnhs&
Q_{\alpha_j,\beta_j}(m_1,n_1)-Q_{\alpha_j,\beta_j}(m_2,n_2) \\
&\dnhs\dnhs&+\ 2(\alpha_j-\alpha)(m_1-m_2)+2(\beta_j-\beta)(n_1-n_2),
\end{eqnarray*}
we have that 
$$
R_2^{(\alpha,\beta)}[-a,a](M_j)\geq R_2^{(\alpha_j,\beta_j)}[0,0](M_j)
$$
when $|\alpha-\alpha_j|< \frac{a}{8(\sqrt{M_j}+1)}$ and
$|\beta-\beta_j|< \frac{a}{8(\sqrt{M_j}+1)}$.
Denote by $B_j\subset\TT^2$ the open set of such $(\alpha,\beta)$.

To summarize,  given any $\lambda>1$, there exists
an $M_j>\lambda$ such that
$$
R_2^{(\alpha,\beta)}[-a,a](M_j)\geq \frac{\log M_j}{\log\log\log M_j}
$$
for all $(\alpha,\beta)\in B_j$. Individually, the sets $B_j$ shrink
to a point as $\lambda\rightarrow\infty$. Note, however, 
that for every fixed $\lambda$ the union
$$
\bigcup_{j: M_j\geq \lambda} B_j
$$
is open and dense in $\TT^2$, and therefore 
$$
B=\bigcap_{\lambda=1}^{\infty} \bigcup_{j: M_j\geq \lambda} B_j
$$
is of second Baire category.

So if $(\alpha,\beta)\in B$, then, given any $\lambda>1$, 
there exists some $M>\lambda$, such that
\vglue12pt
\hfill ${\displaystyle 
R_2^{(\alpha,\beta)}[-a,a](M)\geq \frac{\log M}{\log\log\log M}.
}$ \hfill
\qed
\vglue16pt

Note that the proof remains valid if $\log\log\log$ is replaced by any
slowly increasing positive function $\nu \leq \log\log\log$
with $\nu(M)\rightarrow\infty$ ($M\rightarrow\infty$).

\vglue12pt 9.3. {\it Proof of Theorem}~1.13 (ii).
By virtue of Theorem 1.8, there exists a countable
dense set $\{(\xi_j,\zeta_j)\in\TT^2: j\in\NN \}$ 
for which the pair correlation density of
the forms $O_j:=Q_{\xi_j,\zeta_j}$ is uniform. 
That is, for any $\lambda>1$, we find
some $L_j>\lambda$ such that 
$$
2\pi a -\frac1\lambda < R_2^{(\xi_j,\zeta_j)}[-a,a](L_j)
<2\pi a +\frac1\lambda .
$$
Let $A_j\subset\TT^2$ be the open set of $(\alpha,\beta)$,
which satisfy $|\alpha-\xi_j|< \varepsilon_j$ and
$|\beta-\zeta_j|\break < \varepsilon_j$, where $\varepsilon_j>0$
can be chosen in such a way that
$$
2\pi a -\frac2\lambda < R_2^{(\alpha,\beta)}[-a,a](L_j)
<2\pi a + \frac2\lambda 
$$
for all $(\alpha,\beta)\in A_j$. Now, 
$$
A=\bigcap_{\lambda=1}^{\infty} \bigcup_{j: L_j\geq \lambda} A_j
$$
is again of second category.
We conclude that if $(\alpha,\beta)\in A$, then, given any $\lambda>1$, 
there exists some $L>\lambda$, such that
\vglue12pt 
\hfill${\displaystyle 
2\pi a -\frac2\lambda < R_2^{(\alpha,\beta)}[-a,a](L)
<2\pi a + \frac2\lambda .
}$
\hfill\qed

\vglue17pt 9.4. {\it Conclusion of the proof of Theorem}~1.13. Since $A$ and $B$ are of second Baire category, so is
the intersection $C=A \cap B$.
\hfill \qed

\vglue16pt \centerline{\bf Appendix A. 
 Symmetries}
\vglue12pt 

A.1.  
We have seen in Section \ref{secunipotent}
that in the case when $\alpha,\beta,1$ are linearly
independent over $\QQ$, 
$$
\int_\RR F(u+\i v, 0; \trans(0,0,\alpha,\beta)) h(u) du 
$$
converges for all suitably nice test functions
$F$ on $\Gamma^2\backslash G^2$
to the average of $F$ over $\Gamma^2\backslash G^2$,
as $v\rightarrow 0$.
This is no longer true 
when $\alpha,\beta,1$ are linearly dependent over $\QQ$, i.e., if
we find integers $(k,l,m)\in\ZZ^3-\{(0,0,0)\}$ such that $k\alpha+l\beta+m=0$.
One of $k,l$ must be nonzero, and we will assume in the following 
(without loss of 
generality) that $l\neq 0$, i.e., $\beta=-\frac1l(k\alpha +m)$.

\vglue12pt A.2.  
Suppose $\alpha\notin\QQ$.
For any given function $F\in\C(G^2)$ which is invariant under the left
action of $\Gamma^2$,
we define a function $\tilde F\in\C(G^1)$ by
$$
\tilde F\left(\tau,\phi;\left( \begin{array}{c}
x\\ y\end{array}\right) \right)
=F\left(\tau,\phi;\left( \begin{array}{c}
x\\ -\frac{1}{l}(kx+m) \\ y\\-\frac{1}{l}(ky+m)
\end{array}\right) \right).
$$
Since $\tilde F$ is invariant under the left action of the 
subgroup
$$
\Gamma^1_{2l}=\left\{ (\gamma,\vecn)\in \SL(2,\ZZ)\ltimes\ZZ^2 :
\gamma=\left( \begin{array}{cc} 1 & 0 \\ 0 & 1 \end{array}\right) \bmod  2l,
\;
\vecn=
\vecnull \bmod 2l \right \} \subset \Gamma^1 ,
$$
we can identify $\tilde F$ as a function on $\Gamma^1_{2l}\backslash G^1$.
The congruence group $\Gamma^1_{2l}$ is of finite index in 
$\SL(2,\ZZ)\ltimes\ZZ^2$, and hence
$\Gamma^1_{2l}\backslash G^1$ has finite volume with respect to Haar measure.

\vglue12pt A.3.
If $\alpha=\frac{p}{q}$ and $\beta=\frac{r}{s}$ are rational,
we define instead 
$$
\tilde F(\tau,\phi)=
F(\tau,\phi;\trans(0,0,\tfrac{p}{q},\tfrac{r}{s})) 
$$
which is a function on $G^0$ invariant under the left action
of the subgroup
$$
\Gamma^0_{2qs}
=\left\{ \gamma \in \Gamma_\theta : 
\gamma=\left( \begin{array}{cc} 1 & 0 \\ 0 & 1 \end{array}\right) \bmod  2qs  \right\} .
$$
Again, $\Gamma^0_{2qs}\backslash G^0$ has finite measure.

\vglue12pt A.4. {\it Example} 1.
Consider
the case $\alpha=\beta\notin\QQ$. In order to remove the two-fold
degeneracy we consider the symmetry-reduced set
$$
\{ (m-\alpha)^2 + (n-\alpha)^2 : (m,n)\in\ZZ^2,\; m\geq n \}
$$
whose elements we label by $\lambda_1<\lambda_2<\cdots$.
The pair correlation function of this sequence is now again
Poissonian:

\nonumproclaim{A.5. Theorem}  Assume $\alpha=\beta$ 
is diophantine.
Then
$$
\lim_{\lambda\rightarrow\infty} R_2[a,b](\lambda)  = \frac{\pi}{2} (b-a) .
$$
\endproclaim

Notice that the mean density is now $\frac\pi2$ since we count only
distinct elements.

\vglue12pt A.6. {\it Sketch of the proof}.
The smoothed correlation function is
\begin{eqnarray*}
R_2(\psi_1,\psi_2,h,\lambda)&=&\frac{2}{\pi\lambda}
\sum_{{(m_1,n_1)\in\ZZ^2 \atop m_1\geq n_1}} 
\sum_{{(m_2,n_2)\in\ZZ^2 \atop m_2\geq n_2}}   \\[5pt]
&&\times \
\psi_1\left(\frac{(m_1-\alpha)^2 + (n_1-\alpha)^2}{\lambda} \right) \\[5pt]
&&\times\ 
\psi_2\left(\frac{(m_2-\alpha)^2 + (n_2-\alpha)^2}{\lambda} \right)   \\[5pt]
&&\times\
\hat h((m_1-\alpha)^2 + (n_1-\alpha)^2-(m_2-\alpha)^2 - (n_2-\alpha)^2) .  
\end{eqnarray*}
 
\pagebreak

\noindent 
This is asymptotic for large $\lambda$:
\begin{eqnarray*}
   R_2(\psi_1,\psi_2,h,\lambda)  & \sim&  \frac{1}{2\pi\lambda}
\sum_{{(m_1,n_1)\in\ZZ^2}} 
\sum_{{(m_2,n_2)\in\ZZ^2}}   \psi_1\left(\frac{(m_1-\alpha)^2 + (n_1-\alpha)^2}{\lambda} \right) \\ &&\times \
\psi_2\left(\frac{(m_2-\alpha)^2 + (n_2-\alpha)^2}{\lambda} \right) \\
&&\times \hat h((m_1-\alpha)^2 + (n_1-\alpha)^2-(m_2-\alpha)^2 - (n_2-\alpha)^2) ,
\end{eqnarray*}
since the diagonal terms $m_1=n_1$ or $m_2=n_2$ give  
lower order contributions.
The right-hand side of the above expression is equal to 
$$
\frac{1}{2\pi} \int_\RR 
\Theta_f\left(u+\i \frac{1}{\lambda},0;\trans(0,0,\alpha,\alpha)\right)
\overline{\Theta_g\left(u+\i \frac{1}{\lambda},0;\trans(0,0,\alpha,\alpha)\right)} 
\; h(u)\, du.
$$
The corresponding test function 
$$
\tilde F(\tau,\phi;\trans(x,y))
=\Theta_f(\tau,\phi;\trans(x,x,y,y))
\overline{\Theta_g(\tau,\phi;\trans(x,x,y,y))} 
$$
is a function on $\Gamma^1\backslash G^1$; compare A.2.
Starting from Theorem~7.3 we can apply the same string of arguments
as before. The only main difference is that Lemma~8.2
has to be replaced by the one given below. This yields (compare
the main Theorem 8.1; we assume here  that $\psi_1,\psi_2$ are real-valued)
\begin{eqnarray*}
&&\hskip-.5in \lim_{v\rightarrow 0}
\int_\RR \Theta_f(u+\i v,0;\trans(0,0,y,y))
\overline{\Theta_g(u+\i v,0;\trans(0,0,y,y))} \; h(u)\, du \\
&&= 2\pi \left\{ \pi h(0) + \int_\RR h(u)\, du\right\} 
\int_0^\infty \psi_1(r)\psi_2(r)\, dr , 
\end{eqnarray*}
and hence
$$
\lim_{\lambda\rightarrow\infty}
R_2(\psi_1,\psi_2,h,\lambda) 
= \left\{ \hat h(0) + \frac{\pi}{2} \int_\RR \hat h(s)\,ds \right\}
\int\psi_1(r)\psi_2(r)\, dr ,
$$
as needed.
\hfill\qed

\nonumproclaim{A.7. Lemma}  If $f,g\in\Sw(\RR^2)${\rm ,}
\begin{eqnarray*}
&&\hskip-.5in \frac{1}{\mu(\Gamma^1\backslash G^1)}
\int_{\Gamma^1\backslash G^1} 
\Theta_f(\tau,\phi;\trans(x,x,y,y))
\overline{\Theta_g(\tau,\phi;\trans(x,x,y,y))}\, d\mu \\[4pt]
&&= \iint 
\left\{ f(w_1,w_2)\overline{g(w_1,w_2)}+ f(w_1,w_2)\overline{g(w_2,w_1)}\right\}
\,dw_1\, dw_2 .
\end{eqnarray*}
\endproclaim

When $f(w_1,w_2)=\psi_1(w_1^2+w_2^2)$ and 
$g(w_1,w_2)=\psi_2(w_1^2+w_2^2)$, this yields
$$
\iint\!\!
\left\{ f(w_1,w_2)\overline{g(w_1,w_2)}+ f(w_1,w_2)\overline{g(w_2,w_1)}\right\}\!\!
\,dw_1\, dw_2
= 2\pi \int_0^\infty \psi_1(r)\psi_2(r)\, dr;
$$ 
compare~8.2.

\demo{Proof}
Consider the function
$$
F(\tau,\phi)
= \iint_{\TT^2}
\Theta_f(\tau,\phi;\trans(x,x,y,y))
\overline{\Theta_g(\tau,\phi;\trans(x,x,y,y))}\,dx\,dy .
$$
This function may be viewed as a function on 
$\SL(2,\ZZ)\backslash\SL(2,\RR)$. 
By virtue of Proposition~4.10 one finds that asymptotically
in the cusp ($v\rightarrow\infty$) we have
$$
F(\tau,\phi)
= v^{1/2} \int f_\phi(w,w)\overline{g_\phi(w,w)} dw 
+O_R(v^{-R}) .
$$
It follows from the classical equidistribution of closed horocycles
\cite{Sarnak81}, \cite{Eskin93} in the case of
unbounded test functions (cf.~Proposition 4.3 in \cite{Marklof96})
that as $v\rightarrow 0$
\begin{eqnarray*}
&&\hskip-.25in
\int_0^1 F(u+\i v,0) du \\&&\qquad
= 
\int_0^1 \iint_{\TT^2}
\Theta_f(u+\i v,\phi;\trans(x,x,y,y))
\overline{\Theta_g(\tau,\phi;\trans(x,x,y,y))}
\,dx\,dy\,du 
\end{eqnarray*}
converges to the left-hand side in Lemma~A.7. The right-hand side
of the above equation can, however, be
worked out straightforwardly: The series representation of $\Theta_f$
gives a natural Fourier expansion with respect to $x$ and $u$. The zeroth
Fourier coefficient, which we want to calculate, is given by those 
summands for which
$$
\left\{ \begin{array}{ll}
m_1+n_1=m_2+n_2 & \\
m_1^2+n_1^2=m_2^2+n_2^2 .
\end{array}\right. 
$$
This set of equations is equivalent to
$$
\left\{ \begin{array}{ll}
m_1-m_2=n_2-n_1 & \\
m_1^2-m_2^2=n_2^2-n_1^2 ,
\end{array}\right. 
$$
whose only solutions are obviously $(m_1=m_2,n_1=n_2)$ or
$(m_1=n_2,m_2=n_1)$.
In the limit $v\rightarrow 0$, the zeroth Fourier coefficient
is now easily seen to converge to the right-hand side in Lemma~A.7.
\enddemo

A further special case of interest is the following.

\vglue12pt A.8. {\it Example} 2.
When $\beta=0$ or $\beta=\tfrac12$ we consider the symmetry-reduced
sequences $\lambda_1<\lambda_2<\cdots$ given by the sets
$$
\{ (m-\alpha)^2 + n^2 : (m,n)\in\ZZ^2,\; n\geq 0 \}
$$
or
$$
\{ (m-\alpha)^2 + (n-\tfrac12)^2 : (m,n)\in\ZZ^2,\; n>\tfrac12 \},
$$
respectively. 

\nonumproclaim{A.9. Theorem}  Assume $\alpha$ is diophantine.
Then
$$
\lim_{\lambda\rightarrow\infty} R_2[a,b](\lambda)  = \frac{\pi}{2} (b-a) .
$$
\endproclaim

The proof of this theorem is analogous to that of Theorem~A.5.

\vglue12pt A.10. {\it Example} 3.
If $\alpha=\frac{p}{q}$, $\beta=\frac{r}{s}$ are both rational,
the integral 
$$
\int_\RR \tilde F(u+\i v,0)  h(u) du
$$
of the corresponding test function
$$
\tilde F(\tau,\phi)=
\Theta_f(\tau,\phi;\trans(0,0,\tfrac{p}{q},\tfrac{r}{s})) 
\; \overline{\Theta_g(\tau,\phi;\trans(0,0,\tfrac{p}{q},\tfrac{r}{s}))}
$$
is diverging as $v\rightarrow 0$; one finds in particular that
in this limit
$$
\frac{1}{2qs} \int_0^{2qs} \tilde F(u+\i v,0) du
\sim b_{\alpha,\beta} \log v^{-1}
$$
for some constant $b_{\alpha,\beta}>0$. This follows from
arguments analogous to those given in \cite[Th.~6.1]{Marklof96}. 

Therefore, for $\lambda\rightarrow\infty$,
\begin{eqnarray*}
R_2[0,0](\lambda)&\dnhs=\dnhs &
\frac{1}{\pi \lambda} 
\#\{ (m_1,m_2,n_1,n_2)\in\ZZ^4 : (m_1-\alpha)^2+(n_1-\beta)^2 \leq \lambda,\\
&\dnhs\dnhs&(m_1-\alpha)^2+(n_1-\beta)^2=(m_2-\alpha)^2+(n_2-\beta)^2 \} 
\sim c_{\alpha,\beta} \log\lambda ,
\end{eqnarray*}
for some constant $c_{\alpha,\beta}>0$.
In the case $\alpha=\beta=0$ this yields, of course, Landau's well known result
on the asymptotic number of ways of writing an integer as
a sum of two squares.

\vglue16pt \centerline{\bf Appendix B. Closed connected
subgroups of $\SL(2,\RR)\ltimes\RR^{2k}$}
\vglue12pt

B.1.
Suppose $H$ is a subgroup of $G^k=\SL(2,\RR)\ltimes\RR^{2k}$.
Then 
$$
H=\{ (M;\vecxi) \in G^k: \; M\in L , \; \vecxi\in {\cal C}(M) \}
$$
where $L$ is a subgroup of $\SL(2,\RR)$ and
${\cal C}(M)$ is a family of 
sets, which are suitably chosen such that
$H$ is a group, but are otherwise arbitrary.

\vglue12pt B.2.
Clearly $\Omega={\cal C}(1)$ is a subgroup of $\RR^{2k}$,
because $(1;\vecxi)(1;\vecxi')^{\pm1}=(1;\vecxi\pm \vecxi')$ for all
$\vecxi,\vecxi'\in\Omega$ implies $\vecxi\pm \vecxi'\in\Omega$.

Moreover,
$$
(M;\vecxi)(1;\vecxi')(M;\vecxi)^{-1}=(1;M\vecxi'), 
\hbox{ for any $(M;\vecxi)\in H$,} 
$$
says that if $\vecxi'\in\Omega$, then $M\vecxi'\in\Omega$;
hence $\Omega$ is invariant under the action of $L$.
This means also that $\{1\}\ltimes\Omega$ is a normal subgroup of $H$.
Thus if $(M;\veczeta(M))$ is a 
set of representatives from the coset $(\{1\}\ltimes\Omega)\backslash H$,
$$
(M;\veczeta(M))(M';\veczeta(M'))= (1, \vecsigma(M,M'))(MM';\veczeta(MM'))
$$
with cocycle
$$
\vecsigma(M,M')=\veczeta(M)+M\veczeta(M')-\veczeta(MM') \in\Omega .
$$
We choose $\veczeta(M)$ in such a way that $\veczeta(1)=\vecnull$.

\vglue12pt B.3. 
If $H$ is a closed connected subgroup of $G^k$, then $L$ is
a connected Lie subgroup of $\SL(2,\RR)$. Since all such subgroups
are closed in $\SL(2,\RR)$, $L$ is a closed connected subgroup of
$\SL(2,\RR)$.

\vglue12pt B.4.
Let us assume in the following that the subgroup 
$$
\Psi_0^\RR=
\left(\left( \begin{array}{cc}
1 & \RR \\
0 & 1
\end{array}\right) ; \vecnull \right).
$$
is contained in $H$, and that $L=\SL(2,\RR)$.
Then $$
R=\left( \begin{array}{cc}
0& -1\\
1&  0
\end{array}\right)
\in L 
$$
and thus $(R;\vecxi_0)\in H$ for some vector
$$
\vecxi_0=\left( \begin{array}{c} \vecx_0 \\ \vecy_0 \end{array}\right) .
$$
Since conjugation by
$$
g=\left(1; \tfrac12\left( \begin{array}{c} \vecx_0 - \vecy_0  \\ \vecnull \end{array}\right)\right)
$$
yields
$$
g^{-1} \left(R;\left( \begin{array}{c} \vecx_0 \\ \vecy_0 \end{array}\right)\right) g
=
\left(R;\tfrac12\left( \begin{array}{c} \vecx_0+\vecy_0 \\ \vecx_0+\vecy_0 \end{array}\right)\right)
$$
and
$$
g^{-1}\;\Psi_0^t\; g = \Psi_0^t
$$
for any $t\in\RR$, 
we may assume without loss of generality that $\vecx_0=\vecy_0$
(replace $H$ with $g^{-1}H g$).

Note that
$$
\left(R;\left( \begin{array}{c} \vecx_0 \\ \vecx_0 \end{array}\right)\right)^2
=\left(-1; \left( \begin{array}{c} \vecnull \\ 2\vecx_0 \end{array}\right)\right) 
$$
is in $H$, and so is the conjugate
$$
\left(-1; \left( \begin{array}{c} \vecnull \\ 2\vecx_0 \end{array}\right)\right) 
\left(\left( \begin{array}{cc} 1 & t \\ 0 & 1 \end{array}\right); \vecnull\right)
\left(-1; \left( \begin{array}{c} \vecnull \\ 2\vecx_0 \end{array}\right)\right) 
=
\left(\left( \begin{array}{cc} 1 & t \\ 0 & 1 \end{array}\right); 
\left( \begin{array}{c}  -2t\vecx_0 \\ \vecnull \end{array}\right)\right) .
$$
This implies, however, that 
$$
\left(1,\left( \begin{array}{c} -2t \vecx_0 \\ \vecnull \end{array}\right)\right) \in H
$$
for all $t\in\RR$, and so
$$
\left(1;\left( \begin{array}{c} - \vecx_0 \\ \vecnull \end{array}\right)\right)
\left(R;\left( \begin{array}{c} \vecx_0 \\ \vecx_0 \end{array}\right)\right) 
\left(1;\left( \begin{array}{c} - \vecx_0 \\ \vecnull \end{array}\right)\right)
=(R;\vecnull) \in H.
$$

Because the elements
$$
\left( \begin{array}{cc} 1 & t \\ 0 & 1 \end{array}\right) \; (t\in\RR),
\quad
\left( \begin{array}{cc} 0 & -1 \\ 1 & 0 \end{array}\right)
$$
generate $\SL(2,\RR)$, we find trivially that $\Psi_0^t$ and $(R;\vecnull)$ 
generate
$\SL(2,\RR)\ltimes\{\vecnull\}$, and thus
$H=\SL(2,\RR)\ltimes\Omega$.

\vglue12pt B.5. 
Since $\Omega$ is invariant under the action of $\SL(2,\RR)$
and $H$ is closed and connected,
it is a closed connected subgroup of $\RR^{2k}$, i.e.,
$\Omega$ is a closed  linear subspace. 

\vglue12pt

 B.6.
We conclude that any
closed connected subgroup $H$ of $G^k$, for which $L=\SL(2,\RR)$
and which contains a conjugate of $\Psi_0^\RR$,
is conjugate to $\SL(2,\RR)\ltimes\Omega$, where $\Omega$
is a closed connected subgroup of $\RR^{2k}$. That is, 
$$
H=g_0\; (\SL(2,\RR)\ltimes\Omega) g_0^{-1},
$$
for some $g_0=(M_0;\vecxi_0)\in G^k$.
Because 
$$
(M_0,\vecnull)(\SL(2,\RR)\ltimes\Omega)(M_0,\vecnull)^{-1}
=\SL(2,\RR)\ltimes\Omega
$$
we may take $M_0=1$ without loss of generality, and hence
$$
H= (1;\vecxi_0)(\SL(2,\RR)\ltimes\Omega)(1;-\vecxi_0).
$$

\end{document}